\documentclass[notitlepage]{article}
\usepackage[left=1.3in, right=1.3in, top=1.3in, bottom=1.3in]{geometry}

\usepackage{amsmath,amssymb,amsfonts}
\usepackage{mathtools}
\usepackage{mathptmx}
\usepackage{enumerate}
\usepackage{comment}
\usepackage{color}
\usepackage{lmodern}

\usepackage{graphicx}
\usepackage{multirow}
\usepackage{amsthm}
\usepackage{mathrsfs}
\usepackage{textcomp}
\usepackage{booktabs}
\usepackage{listings}
\DeclareMathOperator*{\argmin}{argmin} 
\usepackage{amsbsy}
\usepackage{hyperref}

\usepackage{algorithmic}
\usepackage[algoruled,boxed,lined]{algorithm2e}

\newtheorem{theorem}{Theorem}[section]
\newtheorem{lemma}[theorem]{Lemma}
\newtheorem{corollary}[theorem]{Corollary}

\theoremstyle{definition}

\newtheorem{remark}{Remark}[section]

\begin{document}

\title{Interior-Point Algorithms for Monotone Linear Complementarity Problem Based on Different Predictor Directions}

\author{Marianna E.-Nagy \and Tibor Ill\'es \and Yurii Nesterov \and Petra Ren\'ata Rig\'o}

\date{\textit{Corvinus Centre for Operations Research at Corvinus Institute for Advanced Studies, Corvinus University of Budapest, Budapest, Hungary}}

\maketitle

\textbf{Abstract.} In this paper, we introduce two parabolic target-space interior-point algorithms for solving monotone linear complementarity problems. The first algorithm is based on a universal tangent direction, which has been recently proposed for linear optimization problems \cite{UTD_LP}.  We prove that this method has the best known worst-case complexity bound. We extend onto LCP its
auto-correcting version \cite{Nesterov_loc_conv}, and prove its local quadratic convergence under a non-degeneracy assumption. In our numerical experiments, we compare the new algorithms with a general method, recently developed for weighted monotone linear complementarity problems \cite{ParLCP}.



\section{Introduction}

Linear complementarity problems (LCPs) have applications in several fields, see \cite{Cottle-Pang-Stone-1992}. It is well-known that linear programming (LP) and linearly constrained (convex) quadratic programming problems are special cases of LCPs. 
The monographs written by Cottle et al. \cite{Cottle-Pang-Stone-1992} and Kojima et al. \cite{kmny1} contain the most important classical results for LCPs. 
General LCPs are NP-complete problems, see \cite{chung1989np}. However, it is known that if the problem's matrix is skew-symmetric (LP) \cite{rtv1,ye1} or positive semidefinite (monotone LCP) \cite{kmy1989}, interior-point algorithms (IPAs) can find an approximate solution of LCPs in polynomial time. 
Starting from a properly selected $\varepsilon$-complementary feasible solution, an exact solution can be obtained using the rounding procedure in strongly polynomial time, which follows for monotone LCPs as a special case of \cite{IT-JP-CR-TT-SIOPT-2000}. Note that for monotone LCPs there are several results about quadratic convergence of the IPAs using Euclidean proximity measures, see \cite{kmny1,anstreicher_ye,potra_Q_conc}. Furthermore, it is also known that some versions of Newton's method have local quadratic convergence for the Nonlinear Complementarity Problem under some non-degeneracy assumptions \cite{DSun}.

In 2008, Nesterov \cite{DAM} presented a new target-following approach for LP, which can start at an arbitrary strictly feasible solution. 
It is based on embedding the problem into a higher-dimensional one, and introducing modeling variables according to the parabolic target space (PTS). The new {\em parabolic barrier function} plays a key role in this algorithm. In \cite{ParLCP}, we generalized this algorithm for weighted monotone LCPs. In \cite{UTD_LP}, we presented a new PTS IPA for LP problems, which is based on a universal tangent direction. In \cite{Nesterov_loc_conv}, it was shown that the algorithm appeared in \cite{UTD_LP} has favorable local behavior under some non-degeneracy assumptions. After that, two new PTS IPAs have been introduced using different search directions in the predictor stage. 

In this paper, we propose two new PTS IPAs for solving monotone LCPs based on the universal tangent direction proposed in \cite{UTD_LP}  and on an auto-correcting version of it based on \cite{Nesterov_loc_conv}. We show that the complexity of the new PTA IPAs coincides with the best-known bounds in the theory of IPAs. 
We present preliminary numerical results, where we compare the two new variants of IPAs with the general method from \cite{ParLCP} for weighted monotone LCPs. 
In the algorithms \cite{ParLCP}, the search directions are obtained by Newton's method as applied to the barrier function, which needs inversion of the Hessians. However, the universal tangent direction and auto-correcting predictor direction used in this paper are less expensive. Another advantage of the auto-correcting IPA is that, although the method operates in a lifted space and it is based on the so-called functional proximity measure, we can prove
its local quadratic convergence under the non-degeneracy assumption. \color{black}

Throughout the paper, we use the following notations.  We use $\mathbb{R}^n_{\oplus}$ and $\mathbb{R}^n_+$ for the positive orthant and its interior.
We denote by $\mathbf{e}$ the vector of all ones of the corresponding size, while $\bar{\mathbf{e}}$ stands for a vector of all ones extended by one additional entry.
In general, the boldface small letters denote finite-dimensional vectors, while the real numbers and coordinates of vectors are denoted by small letters. All arithmetic operations and relations involving vectors, like $\sqrt{\mathbf{x}}$, $\mathbf{x} \mathbf{s}$, $\mathbf{x} / \mathbf{s}$ for $\mathbf{x}, \mathbf{s} \in \mathbb{R}^n$, are understood in the component-wise sense. The scalar products and the norms are defined in the standard way:
$\mathbf{x}^T \mathbf{s} = \sum\limits_{i=1}^n x_i \, s_i, \; \| \mathbf{x} \|^2 = \mathbf{x}^T \mathbf{x}, \; \| \mathbf{x} \|_1 =\sum\limits_{i=1}^n |x_i|, \; \| \mathbf{x} \|_{\infty} =\max\limits_{1\leq i \leq n}  |x_i|, \mathbf{x}, \mathbf{s} \in \mathbb{R}^n.$
Note that  $X$ denotes the diagonal matrix containing the components of the vector $\mathbf{x}$ on the main diagonal. In our analysis, we often use the following functions \cite{NN,YN-LN-2018}:
\begin{equation}\label{omegas}
\omega(t) = t - \ln (1+t): \; \mathbb{R}_+ \to \mathbb{R}, \qquad \omega_*(t) = - t - \ln (1-t): \; [0,1) \to \mathbb{R}.
\end{equation}

\section{Parabolic target space}

Consider the LCPs 
 {\renewcommand{\theequation}{$LCP$}
 \begin{equation}
-M \mathbf{x} + \mathbf{s} = \mathbf{q}, \quad \mathbf{x}, \mathbf{s} \geq \mathbf{0}, \quad
\mathbf{x} \, \mathbf{s} = \mathbf{0}, \label{LCP}
\end{equation}
\addtocounter{equation}{-1}}where $M\in \mathbb{R}^{n\times n}$ is a given matrix, $\mathbf{q} \in \mathbb{R}^{n}$ is a given vector and $n$ is a natural number. \medskip

We assume that the matrix $M$ is positive semidefinite. Denote by
\[
{\mathcal{F}} = \{(\mathbf{x}, \mathbf{s}) \in \mathbb{R}_{\oplus}^{n} \times \mathbb{R}_{\oplus}^{n}: -M \mathbf{x} + \mathbf{s} = \mathbf{q} \}, \quad
\mathcal{F}^{+} = \{(\mathbf{x}, \mathbf{s}) \in \mathbb{R}_{+}^{n} \times \mathbb{R}_{+}^{n}: -M \mathbf{x} + \mathbf{s} = \mathbf{q}\}\]
and  ${\mathcal{F}^{*}} = \{(\mathbf{x}, \mathbf{s}) \in \mathcal{F}:\mathbf{x} \, \mathbf{s} = {\mathbf{0}} \}$ the sets of feasible,
strictly feasible, and optimal solutions to LCP. Our main assumption is ${\cal F}^+ \neq \emptyset$.

Following \cite{ParLCP}, we relax the nonlinear equation of (\ref{LCP}) and define the following {\textsl convex feasibility problem }
{\renewcommand{\theequation}{$CFP$}
\begin{equation}\label{def-CFP}
  -M \mathbf{x} + \mathbf{s} = \mathbf{q}, \quad \mathbf{x}, \mathbf{s} \ge \mathbf{0}, \quad \mathbf{x} \, \mathbf{s} \ge \mathbf{v}^2 \quad \hbox{and} \quad v_0 \ge \mathbf{x}^T \mathbf{s}.
\end{equation}\addtocounter{equation}{-1}}

When $(\mathbf{x},\mathbf{s})$ is a solution of (\ref{def-CFP}), then $v_0 \ge \| \mathbf{v} \|^2$. Therefore, as in \cite{DAM,ParLCP,UTD_LP}, we introduce the variables $\mathbf{w} = (v_0, \mathbf{v}) \in \mathbb{R}^{1+n}$ from the {\textsl parabolic target space (PTS)} 
\[ {\mathcal{F}_w} = \{ \mathbf{w} = (v_0, \mathbf{v}) \in \mathbb{R}^{1+n}: v_0 \ge \|\mathbf{v}\|^2 \}.
\]
Clearly, $\mathcal{F}_{w}$ is a convex set. We will use notation $\mathcal{F}_{w}^+$ to denote its interior.

The main idea of our approach is taking a sequence from $\mathcal{F}_w^+$ which tends to the origin, and solving approximately the corresponding (\ref{def-CFP}) problems.
Let us denote the set of all vectors satisfying the system of convex inequalities (\ref{def-CFP}) by
\[
\mathcal{F}_{z} = \left\{ \mathbf{z} = (\mathbf{x},\mathbf{s},\mathbf{w}) \in \mathcal{F} \times \mathcal{F}_{w} : \;  \mathbf{x} \, \mathbf{s} \geq \mathbf{v}^2, \;   \mathbf{x}^T  \mathbf{s} \leq v_0 \right\}.
\]
Since $M$ is a positive semidefinite matrix, the restriction of function $\mathbf{x}^T \mathbf{s}$ onto the set $\mathcal{F}$ is convex in $\mathbf{x}$:
\[
\mathbf{x}^T \mathbf{s} = \mathbf{x}^T(M\mathbf{x}+\mathbf{q}) = \mathbf{x}^TM\mathbf{x} + \mathbf{q}^T\mathbf{x}.
\] 
Therefore, the set $\mathcal{F}_z$ is convex.

From $\mathcal{F}^+ \neq \emptyset$ follows that the relative interior of $\mathcal{F}_{z}$, denoted by $\mathcal{F}_{z}^+$, is nonempty. It admits a standard self-concordant barrier 
\begin{equation}\label{sc_barrier}
F(\mathbf{z}) \, = \, F(\mathbf{x},\mathbf{s}, \mathbf{w}) = - \ln \left(v_0 - \mathbf{x}^T \mathbf{s} \right)-
\sum_{i=1}^n \, \ln \left( x_i s_i  - v^2_i \right), \quad \mathbf{z} \in \mathcal{F}_z^+
\end{equation}
with barrier parameter $\nu_F = 2n+1$. 

Let us define the following residuals:
\begin{equation}\label{res}
r_0(\mathbf{z})  =  v_0 - \mathbf{x}^T \mathbf{s}, \quad r_i(\mathbf{z}) \; = \; x_i \, s_i  - v^2_i, \quad i = 1, \dots, n.
\end{equation}
Hence, $\mathbf{r}(\mathbf{z}) \in \mathbb{R}^{1+n}_{\oplus}$. 
Note that
\begin{equation}\label{res_av}
\frac{\bar{\mathbf{e}}^T \mathbf{r}(\mathbf{z})}{n+1}= \frac{v_0-\|\mathbf{v}\|^2}{n+1} =:\rho(\mathbf{w}). 
\end{equation}
With notation (\ref{res}), the self-concordant barrier (\ref{sc_barrier}) can be written as
\begin{equation}\label{sc_barrier2}
F(\mathbf{z})  =  - \sum\limits_{i=0}^n \ln (r_i(\mathbf{z})).
\end{equation}

Let us define now the control barrier function
\begin{equation}\label{def-Phi}
  \phi(\mathbf{w}) = \min\limits_{(\mathbf{x},\mathbf{s}): (\mathbf{x},\mathbf{s},\mathbf{w}) \in \mathcal{F}_z} F(\mathbf{x},\mathbf{s},\mathbf{w}) \; : \quad \mathbb{R}^{1+n} \to \mathbb{R},
\end{equation}
with notation  $\mathbf{z}(\mathbf{w}) = \left( \mathbf{x}(\mathbf{w}),\mathbf{s}(\mathbf{w}),\mathbf{w} \right)$  for the optimal solution of (\ref{def-Phi}). As in \cite{DAM,ParLCP}, it has a closed form representation. 

\begin{theorem}[Theorem 3.2 in \cite{ParLCP}] \label{th-Phi} Let $M$ be positive semidefinite and $\mathcal{F}^+ \neq \emptyset$. Then ${\rm dom} \, \phi = \mathcal{F}_{w}^+ \neq \emptyset$ and the optimization problem (\ref{def-Phi}) has a unique solution. The optimal vector $(\mathbf{x}(\mathbf{w}),\mathbf{s}(\mathbf{w}))$ is the unique solution of the following system: 
\begin{equation}\label{eq-XOpt}
\begin{array}{rl}
-M \mathbf{x} + \mathbf{s} &= \mathbf{q}, \quad \mathbf{x},\mathbf{s} \geq \mathbf{0},\\ 
    \mathbf{x} \,  \mathbf{s} &= \mathbf{v}^2 + \frac{v_0 - \| \mathbf{v} \|^2}{n+1} \, \mathbf{e}.
    \end{array}
\end{equation}

Moreover, for all $\mathbf{w} \in \mathcal{F}_w^+$, we have
  \begin{equation*}
    \phi(v_0, \mathbf{v}) = -(n+1) \, \ln \frac{v_0 - \| \mathbf{v} \|^2}{n+1}.
  \end{equation*}
\end{theorem}
\begin{corollary}
Under assumptions of Theorem \ref{th-Phi}, the set ${\cal F}^*$ is nonempty and bounded.
\end{corollary}
\proof{Proof.}
Indeed, let $(\hat{\mathbf{x}},\hat{\mathbf{s}}) \in {\cal F}^+$. For for any $z = (\mathbf{x},\mathbf{s}, (v_0, \mathbf{v})) \in {\cal F}_z$, we have
\[
0 \leq ( \hat{\mathbf{x}} - \mathbf{x})^T(\hat{\mathbf{s}} - \mathbf{s}) \leq  \hat{\mathbf{x}}^T \hat{\mathbf{s}} - \mathbf{x}^T 
\hat{\mathbf{s}} -  \hat{\mathbf{x}}^T \mathbf{s} + v_0.
\]
Therefore, for any sequence $\mathbf{w}_k \in {\cal F}_w^+$ converging to the origin, the sequence $(\mathbf{x}(\mathbf{w_k}),\mathbf{s}(\mathbf{w_k}))$ is bounded and its limit points belong to ${\cal F}^*$.
 \endproof

Let $\hat{\mathbf{r}}(\mathbf{z}) \in \mathbb{R}_+^{1+n}$ be defined as
\begin{equation}\label{t}
\hat{\mathbf{r}}(\mathbf{z}) = \sqrt{\frac{\mathbf{r}(\mathbf{z})}{\rho(\mathbf{w})}}.
\end{equation}
We have
\begin{equation}\label{normt}
\bar{\mathbf{e}}^T \hat{\mathbf{r}}^2(\mathbf{z}) = \| \hat{ \mathbf{r}}(\mathbf{z})\|^2 = n+1.
\end{equation}

By definition $\hat{\mathbf{r}}(\mathbf{z}) = \bar{\mathbf{e}}$ if and only if $(\mathbf{x},\mathbf{s})=(\mathbf{x}(\mathbf{w}),\mathbf{s}(\mathbf{w}))$.
Note that we cannot compute $(\mathbf{x}(\mathbf{w}),\mathbf{s}(\mathbf{w}))$ in a closed form.  However, using (\ref{t}) we can define the following measures of closeness of point $(\mathbf{x},\mathbf{s})$ to $(\mathbf{x}(\mathbf{w}),\mathbf{s}(\mathbf{w}))$:
\begin{equation}\label{def-Rho}
\zeta_0^2(\mathbf{z})\; = \; \left\| \hat{\mathbf{r}}(\mathbf{z}) - \frac{1}{\hat{\mathbf{r}}(\mathbf{z})} \right\|^2, \quad
 \zeta_1(\mathbf{z})=\left\|\frac{1}{\hat{\mathbf{r}}^2(\mathbf{z})}-\bar{\mathbf{e}} \right\|, \quad \zeta_2(\mathbf{z}) = \| \hat{ \mathbf{r}}^2(\mathbf{z})-\mathbf{\bar{e}} \|,
\end{equation}
\begin{equation}\label{def-delta}
\delta(\mathbf{z}) \; = \; \frac{\zeta_0^2(\mathbf{z})}{\zeta_1(\mathbf{z})} \quad (\mbox{for $\zeta_1(\mathbf{z}) = 0$, we define $\delta(\mathbf{z}) = 0$}).
\end{equation}

Using $\delta(\mathbf{z})$, we can define the following neighbourhood with parameter $\beta > 0$:
\begin{equation*}
\mathcal{N}_{\delta}(\beta)= \{\mathbf{z} \in \mathcal{F}_z^+: \delta(\mathbf{z}) \leq \beta \}.
\end{equation*}
Another possibility to measure the closeness of $(\mathbf{x},\mathbf{s})$ to $(\mathbf{x}(\mathbf{w}),\mathbf{s}(\mathbf{w}))$ is based on the definition of the control barrier function $\phi$.
For the unique optimal solution $(\mathbf{x}(\mathbf{w}),\mathbf{s}(\mathbf{w}))$ of problem (\ref{def-Phi}), we have $F(\mathbf{x}(\mathbf{w}),\mathbf{s}(\mathbf{w}),\mathbf{w}) = \phi(\mathbf{w})$ and $F(\mathbf{z})=F(\mathbf{x},\mathbf{s},\mathbf{w}) \geq \phi(\mathbf{w})$ for all $\mathbf{z}  \in \mathcal{F}_z^+$. Therefore, we can define the {\em functional proximity measure}
\begin{equation}\label{def-FProx}
\Psi(\mathbf{z}) = F(\mathbf{z}) - \phi(\mathbf{w}) =   - \sum\limits_{i=0}^n \ln  r_i(\mathbf{z})   + (n+1) \ln \rho(\mathbf{w}) = - \sum\limits_{i=0}^n \ln \hat{r}_i^2(\mathbf{z})  \geq 0.
\end{equation}
Using this measure, we can introduce the following neighbourhood:
\begin{equation*}
\mathcal{N}_{\Psi}(\tau) = \{\mathbf{z}  \in \mathcal{F}_z^+: \Psi(\mathbf{z}) \leq \tau\}, \quad \tau \geq 0.
\end{equation*}

\section{Interior-point algorithm based on universal tangent direction}

If $(\mathbf{x}, \mathbf{s}) \in \mathcal{F}^+$, then we can always find a vector $\mathbf{w}(\mathbf{x},\mathbf{s}) \in \mathcal{F}_w^+ $, such that $(\mathbf{x},\mathbf{s}) = (\mathbf{x}(\mathbf{w}(\mathbf{x},\mathbf{s})), \mathbf{s}(\mathbf{w}(\mathbf{x},\mathbf{s})))$. Indeed, define
$
\xi(\mathbf{x},\mathbf{s}) = \min\limits_{1 \leq i \leq n} \, x_i \, s_i \, > \, 0
$
and
\begin{equation}\label{def-W(x)}
\mathbf{v}(\mathbf{x},\mathbf{s}) = \Big[\mathbf{x} \,\mathbf{s} - \xi(\mathbf{x},\mathbf{s}) \, \mathbf{e} \Big]^{1/2}, \qquad
v_0(\mathbf{x},\mathbf{s}) = \mathbf{x}^T \mathbf{s} + \xi(\mathbf{x},\mathbf{s}).
\end{equation}
In this way, we have 
\begin{equation*}
  \xi(\mathbf{x},\mathbf{s}) = \frac{v_0(\mathbf{x},\mathbf{s}) - \| \mathbf{v}(\mathbf{x},\mathbf{s}) \|^2}{n+1}.
\end{equation*}
Then for $\mathbf{w}(\mathbf{x},\mathbf{s}) = (v_0(\mathbf{x},\mathbf{s}), \mathbf{v}(\mathbf{x},\mathbf{s}))$, we have $(\mathbf{x},\mathbf{s}) = (\mathbf{x}(\mathbf{w}(\mathbf{x},\mathbf{s})), \mathbf{s}(\mathbf{w}(\mathbf{x},\mathbf{s}))).$

Let us define the search directions by the following system of linear equations:
\begin{equation}\label{ASD}
\begin{array}{rl}
-M \Delta \mathbf{x} + \Delta \mathbf{s} &= \mathbf{0} \\
 S \Delta \mathbf{x} + X \Delta \mathbf{s}&= \mathbf{a},
\end{array}
\end{equation}
where $X = \mathrm{diag}(\mathbf{x})$ and $S= \mathrm{diag}(\mathbf{s})$ and $\mathbf{a}\in\mathbb{R}^n$ is dependent on particular strategy. 

In this paper, we reconsider directions proposed in \cite{UTD_LP,Nesterov_loc_conv} for LP problems. 
Since bilinear constraints of system (\ref{eq-XOpt}) are the same in the LP and LCP cases, we have the same right-hand side in system (\ref{ASD}) for defining the search directions. 

Let us look at the following search directions.
\begin{itemize}
\item The {\sc corrector} search direction was used in \cite{UTD_LP} for LP problems. In this case,
\begin{equation}\label{corr}
\mathbf{a}_{c}(\mathbf{z}) = \rho(\mathbf{w}) \mathbf{e} - \mathbf{x}\mathbf{s} + \mathbf{v}^2.
\end{equation}

\item The {\sc universal tangent} direction was  used in \cite{UTD_LP} for LP problems. In this case
\begin{equation}\label{UTD}
\mathbf{a}_{ut}(\mathbf{z}) = \left(\frac{\| \mathbf{v} \|^2}{n+1} - \rho(\mathbf{w}) \right)\mathbf{e} - 2 \mathbf{v}^2 = \frac{v_0}{n+1} \mathbf{e} - 2\mathbf{x}(\mathbf{w})\mathbf{s}(\mathbf{w}).
\end{equation}

\item The {\sc auto-corrector} direction was proposed in \cite{Nesterov_loc_conv} for LP problems. In this case
\begin{equation}\label{autocorr1}
\mathbf{a}_{ac}(\mathbf{z}) = \mathbf{a}_{ut}(\mathbf{z})+\mathbf{a}_{c}(\mathbf{z})= \frac{\|\mathbf{v}\|^2}{n+1} \mathbf{e} - \mathbf{v}^2 - \mathbf{x}\mathbf{s}.
\end{equation}
\end{itemize}

In this paper, we consider corrector search direction $\mathbf{a}_{c}(\mathbf{z})$ defined by (\ref{corr}). For the predictor step, we analyze two special cases:
\begin{itemize}
\item for $\mathbf{a}_{p}(\mathbf{z}) = \mathbf{a}_{ut}(\mathbf{z})$ from (\ref{UTD}), we call the algorithm Universal Tangent Direction-based IPA with shortcut UTD IPA. 
\item for $\mathbf{a}_{p}(\mathbf{z}) = \mathbf{a}_{ac}(\mathbf{z})$ from (\ref{autocorr1}), we call the corresponding algorithm Auto-Corrector IPA with shortcut AC IPA.
\end{itemize}

Below, we present our new {\textsl Parabolic Target-Space IPA} for monotone LCPs. It performs a predictor and a sequence of corrector steps with $\mathbf{w}$-component of points $\mathbf{z}$ being unchanged. 

\begin{algorithm}
\caption{IPA for monotone LCPs}\label{alg:PTS}
\begin{algorithmic}
\REQUIRE
 the accuracy parameter $\varepsilon > 0$, $\beta \in \left(0,{1 \over 3}\right)$, $\tau = \omega_*\left({2\beta \over 1 - \beta}\right)$,\\ 
 \hspace{8mm} the initial point $\mathbf{z}^{(0)} = (\mathbf{x}^{(0)}, \mathbf{s}^{(0)}, \mathbf{w}^{(0)}) \in   \mathcal{N}_{\delta}(\beta)$. \\
 \STATE $\mathbf{x}:=\mathbf{x}^{(0)}; 
 \quad \mathbf{s}:= \mathbf{s}^{(0)}; \quad \mathbf{w}:=\mathbf{w}^{(0)}$.
\WHILE{$v_0  > \varepsilon$}
\STATE \textbf{Predictor step:}
\STATE \quad  Compute $(\Delta \mathbf{x}, \Delta \mathbf{s})$ satisfying (\ref{ASD}) with $\mathbf{a}_{p}(\mathbf{z})$.
\STATE \quad Compute $\alpha_{p} := \displaystyle\max\limits_{\alpha\in[0,1]}\left\{\alpha : \left(\mathbf{x}+\alpha \Delta \mathbf{x}, \mathbf{s}+\alpha \Delta \mathbf{s},(1-\alpha)  \mathbf{w} \right) \in \mathcal{N}_{\Psi}(\tau)\right\}$. 
\STATE \quad $\mathbf{x}:= \mathbf{x}+\alpha_{p} \Delta \mathbf{x}; \quad \mathbf{s}:= \mathbf{s}+\alpha_{p} \Delta \mathbf{s}; \quad \mathbf{w}:=(1-\alpha_{p})\mathbf{w}.$
\WHILE{$\mathbf{z} \notin \mathcal{N}_{\delta}(\beta)$}
\STATE \textbf{Corrector step:}
\STATE \quad Compute $(\Delta \mathbf{x}, \Delta \mathbf{s})$ satisfying (\ref{ASD}) with $\mathbf{a}_{c}(\mathbf{z})$. 
\STATE \quad $\alpha_{c} := \displaystyle \argmin
 \limits_{\alpha \in [0,1]} \{F(\mathbf{x}+\alpha \Delta \mathbf{x}, \mathbf{s}+\alpha \Delta \mathbf{s}, \mathbf{w}): (\mathbf{x}+\alpha \Delta \mathbf{x},\mathbf{s}+\alpha \Delta \mathbf{s}) \in \mathcal{F}^{+}\}$.
\STATE \quad $\mathbf{x}:=\mathbf{x}+\alpha_{c} \Delta \mathbf{x}; \quad \mathbf{s}:=\mathbf{s}+\alpha_{c} \Delta \mathbf{s}.$
\ENDWHILE
\ENDWHILE
\end{algorithmic}
\end{algorithm}

\section{Complexity analysis}\label{sc-Complex}

In this section, we follow the main ideas of \cite{UTD_LP} for LP with slight modifications required by LCP.
Our search directions are defined by the system (\ref{ASD}). 

The following estimates are well known in the theory of LCPs, however, we provide the proof for the reader's convenience.

\begin{lemma}\label{search-dir}
Let $\Delta \mathbf{x}$ and $\Delta \mathbf{s}$ be the solution of system (\ref{ASD}), where $M$ is a positive semidefinite matrix. Then, we have
\begin{equation}\label{search-dir1}
\| \Delta \mathbf{x} \Delta \mathbf{s}\|_{\infty}\leq  \frac{1}{4} \left\|\frac{\mathbf{a}}{\sqrt{\mathbf{x}\mathbf{s}}} \right\|^2, \qquad \| \Delta \mathbf{x} \Delta \mathbf{s}\|_1 \leq  \frac{1}{2} \left\|\frac{\mathbf{a}}{\sqrt{\mathbf{x}\mathbf{s}}} \right\|^2,
\end{equation}
\begin{equation}\label{search-dir2}
\| \Delta \mathbf{x} \Delta \mathbf{s}\| \leq 
\frac{1}{2\sqrt{2}} \left\|\frac{\mathbf{a}}{\sqrt{\mathbf{x}\mathbf{s}}} \right\|^2
\end{equation}
and
\begin{equation}\label{search-dir3}
\Delta \mathbf{x}^T \Delta \mathbf{s} \leq \frac{1}{4} \left\|\frac{\mathbf{a}}{\sqrt{\mathbf{x}\mathbf{s}}} \right\|^2.
\end{equation}
\end{lemma}
\proof{Proof.}
Since $M$ is positive semidefinite, we have
$\Delta\mathbf{x}^T\Delta\mathbf{s}= \Delta\mathbf{x}^T M\Delta\mathbf{x} \geq 0.$
With the notation $\mathcal{I}_+=\{i:\Delta x_i\Delta s_i>0\}$ and $\mathcal{I}_-=\{i:\Delta x_i\Delta s_i<0\}$, we have
\begin{equation}\label{ineq:I-}
\sum_{i\in\mathcal{I}_-} |\Delta x_i\Delta s_i|\le  \sum_{i\in\mathcal{I}_+} \Delta x_i\Delta s_i \le \frac{1}{4} \left\|\sqrt{\frac{\mathbf{x}}{\mathbf{s}}}\Delta\mathbf{x} + \sqrt{\frac{\mathbf{s}}{\mathbf{x}}}\Delta\mathbf{s}\right\|^2=\frac{1}{4}\left\|\frac{\mathbf{a}}{\sqrt{\mathbf{x}\mathbf{s}}}\right\|^2,
\end{equation}
where the last equation comes from rescaling the second equation of the system (\ref{ASD}).

This immediately proves the first inequality of the lemma. We obtain the second inequality from the identity $\|\Delta\mathbf{x}\Delta\mathbf{s}\|_1 =\sum_{i\in\mathcal{I}_-} |\Delta x_i\Delta s_i|+  \sum_{i\in\mathcal{I}_+} \Delta x_i\Delta s_i$.
For the last inequality, we use the relation $\Delta\mathbf{x}^T\Delta\mathbf{s}\le\sum_{i\in\mathcal{I}_+} \Delta x_i\Delta s_i$.

Finally,
$$\left\| \Delta\mathbf{x}\Delta\mathbf{s} \right\|^2
=\sum_{i\in\mathcal{I}_-} \left(\Delta x_i\Delta s_i\right)^2+  \sum_{i\in\mathcal{I}_+} \left(\Delta x_i\Delta s_i\right)^2
\le \left(\sum_{i\in\mathcal{I}_-} |\Delta x_i\Delta s_i|\right)^2+  \left(\sum_{i\in\mathcal{I}_+} \Delta x_i\Delta s_i\right)^2,$$
which together with (\ref{ineq:I-}) proves (\ref{search-dir3}).
 \endproof

\subsection{Corrector stage}

Let us bound the number of corrector steps inside the main iteration. The proof is similar to Lemma 5.1 in \cite{UTD_LP}.

\begin{lemma}\label{nr_corr_steps}
Let $\mathbf{z} \in \mathcal{N}_{\psi}(\tau),$  where $\tau > 0$. Then, each main iteration needs at most
\begin{equation}\label{eq-NUp}
 \left\lfloor  \frac{\tau}{\omega\left(\frac{\beta}{\sqrt{1+2\beta}} \right)}    \right\rfloor + 1
\end{equation}
corrector steps in order to have $\mathbf{z} \in \mathcal{N}_{\delta}(\beta)$, where $\omega$ is defined in (\ref{omegas}).  
\end{lemma}
\proof{Proof.}
After the corrector step, we need to get into the neighbourhood $\mathcal{N}_{\delta}(\beta)$ with $\delta(\mathbf{z}) \leq \beta$. Let us give a lower bound on the decrease of $\delta(\mathbf{z})$ after one corrector step. Let us analyze an effect of
a single damped Newton step applied to the function $f(\alpha)= F(\mathbf{x}+\alpha \Delta \mathbf{x},\mathbf{s} + \alpha \Delta \mathbf{s}, \mathbf{w}) = F(\mathbf{z}+\alpha \Delta \mathbf{z})$, where $\Delta \mathbf{z} =(\Delta \mathbf{x}, \Delta \mathbf{s}, \mathbf{0})$. Using (\ref{ASD}) with (\ref{corr}) and the definitions of $\mathbf{r}(\mathbf{z})$ in (\ref{res}) and $\rho(\mathbf{w})$ in (\ref{res_av}), we obtain 
\begin{align*}
r_0(\mathbf{z}+\alpha \Delta \mathbf{z}) &= (1-\alpha) r_0(\mathbf{z}) + \alpha \rho(\mathbf{w})  - \alpha^2 \Delta \mathbf{x}^T \Delta \mathbf{s},\\
r_i(\mathbf{z}+\alpha \Delta \mathbf{z}) &= (1-\alpha) r_i(\mathbf{z}) + \alpha \rho(\mathbf{w}) + \alpha^2 \Delta x_i \Delta s_i,\quad i=1,\ldots,n.
\end{align*}
Using (\ref{sc_barrier2}) and (\ref{t}), we  have 
\begin{align*}\label{faplha}
f(\alpha)  &= - (n+1) \ln \rho(\mathbf{w}) - \ln  \left[ \hat r_0^2(\mathbf{z}) + \alpha \left(1-\hat r_0^2(\mathbf{z})\right)  -  \alpha^2 \frac{\Delta \mathbf{x}^T \Delta \mathbf{s}}{\rho(\mathbf{w})}\right]\\
 &-\sum_{i=1}^{n} \ln \left[ \hat r_i^2(\mathbf{z}) + \alpha  \left(1-\hat r_i^2(\mathbf{z})\right) + \frac{\alpha^2 \Delta x_i \Delta s_i}{\rho(\mathbf{w})} \right]
\end{align*}
and
\begin{equation*}
f'(\alpha) = -\frac{ 1-\hat r_0^2(\mathbf{z}) - 2 \alpha \frac{\Delta \mathbf{x}^T \Delta \mathbf{s}}{\rho(\mathbf{w})}}{ \hat r_0^2(\mathbf{z}) + \alpha  \left(1-\hat r_0^2(\mathbf{z})\right) - \alpha^2 \frac{\Delta \mathbf{x}^T \Delta \mathbf{s}}{\rho(\mathbf{w})} } - \sum_{i=1}^n \frac{1-\hat r_i^2(\mathbf{z}) + 2 \alpha \frac{\Delta x_i \Delta s_i}{\rho(\mathbf{w})}}{ \hat r_i^2(\mathbf{z}) + \alpha  \left(1-\hat r_i^2(\mathbf{z})\right) + \frac{\alpha^2 \Delta x_i \Delta s_i}{\rho(\mathbf{w})}}.
\end{equation*}
Then, using (\ref{normt}) and (\ref{def-Rho}), we have 
\begin{equation}\label{falphader0}
f'(0) =w
- \sum_{i=0}^{n} \frac{1-\hat r_i^2(\mathbf{z})}{ \hat r_i^2(\mathbf{z})}  =  \mathbf{e}^T\left(\mathbf{e}-\frac{1}{\hat r^2(\mathbf{z})}\right) = -\zeta_0^2(\mathbf{z}).
\end{equation}
Furthermore, 
\begin{align*}
f''(\alpha) &= \frac{\frac{2\Delta \mathbf{x}^T \Delta \mathbf{s}}{\rho(\mathbf{w})}}{ \hat r_0^2(\mathbf{z}) + \alpha  \left(1-\hat r_0^2(\mathbf{z})\right) -\alpha^2 \frac{\Delta \mathbf{x}^T \Delta \mathbf{s}}{\rho(\mathbf{w})} } + \left(\frac{ 1-\hat r_0^2(\mathbf{z}) -  2 \alpha \frac{\Delta \mathbf{x}^T \Delta \mathbf{s}}{\rho(\mathbf{w})}}{ \hat r_0^2(\mathbf{z}) + \alpha  \left(1-\hat r_0^2(\mathbf{z})\right) -  \alpha^2 \frac{\Delta \mathbf{x}^T \Delta \mathbf{s}}{\rho(\mathbf{w})} }
\right)^2 \\
&-\sum_{i=1}^n \frac{\frac{2 \Delta x_i \Delta s_i}{\rho(\mathbf{w})}}{\hat r_i^2(\mathbf{z}) + \alpha  \left(1-\hat r_i^2(\mathbf{z})\right) + \frac{\alpha^2 \Delta x_i \Delta s_i}{\rho(\mathbf{w})}}
+  \sum_{i=1}^n \left(\frac{
1-\hat r_i^2(\mathbf{z}) + 2 \alpha \frac{\Delta x_i \Delta s_i}{\rho(\mathbf{w})}}{ \hat r_i^2(\mathbf{z}) + \alpha  \left(1-\hat r_i^2(\mathbf{z})\right) + \frac{\alpha^2 \Delta x_i \Delta s_i}{\rho(\mathbf{w})}}\right)^2
\end{align*}
and from (\ref{t}) and (\ref{def-Rho}), we have
\begin{align}\label{falphader20}
f''(0) &= \sum_{i=0}^n \left(\frac{1}{\hat r_i^2(\mathbf{z})} - 1\right)^2 -\sum_{i=1}^n \frac{2 \Delta x_i \Delta s_i}{\rho(\mathbf{w}) \hat r_i^2(\mathbf{z})}  + \frac{2 \Delta \mathbf{x}^T \Delta  \mathbf{s}}{\rho(\mathbf{w}) \hat r_0^2(\mathbf{z})} \nonumber\\
&=  \zeta_1^2(\mathbf{z}) + \frac{2}{\rho(\mathbf{w})}\sum_{i=1}^n \Delta x_i \Delta s_i \left(\frac{1}{\hat r_0^2(\mathbf{z})} - \frac{1}{\hat r_i^2(\mathbf{z})} \right) \nonumber\\
&\leq  \zeta_1^2(\mathbf{z}) + \frac{2} {\rho(\mathbf{w})}\sum_{i=1}^n |\Delta x_i \Delta s_i|   \left( \left| \frac{1}{\hat r_0^2(\mathbf{z})}-1 \right| +\left| \frac{1}{\hat r_i^2(\mathbf{z})} -1 \right| \right) \nonumber\\
&\leq  
\zeta_1^2(\mathbf{z}) +  \frac{4 \zeta_1(\mathbf{z})}{\rho(\mathbf{w})} \sum_{i=1}^n |\Delta x_i \Delta s_i|.  
\end{align}

Using (\ref{search-dir1}) of Lemma \ref{search-dir} with $\mathbf{a} = \mathbf{a}_{c}(\mathbf{z})$ given in  (\ref{corr}), we have
\begin{equation}\label{DeltaxDeltas2}
{2 \over \rho(\mathbf{w})} \sum\limits_{i=1}^n \left| \Delta x_i \Delta s_i \right|  \leq  
  \sum\limits_{i=1}^n {( r_i(\mathbf{z}) - \rho(\mathbf{w}))^2 \over \rho(\mathbf{w}) x_is_i} \leq  \sum\limits_{i=0}^n {( r_i(\mathbf{z}) - \rho(\mathbf{w}))^2 \over \rho(\mathbf{w}) r_i(\mathbf{z})} \;=   \;  \zeta_0^2(\mathbf{z}).
\end{equation}

From (\ref{falphader20}) and (\ref{DeltaxDeltas2}), we derive the following upper bound:
\begin{equation*}
f''(0) \leq \zeta_1^2(\mathbf{z}) + 2 \zeta_0^2(\mathbf{z}) \zeta_1(\mathbf{z}).
\end{equation*}

Thus, we proved following lower bound for the Newton decrement of function  $f(\cdot)$:
\begin{equation*}
\lambda  :=  \sqrt{ {(f'(0))^2 \over f''(0)}} \; \geq  \; {\zeta_0^2(\mathbf{z}) \over \sqrt{\zeta_1^2(\mathbf{z}) + 
 2 \zeta_0^2(\mathbf{z}) \zeta_1(\mathbf{z})}} \; = \; {\delta(\mathbf{z}) \over \sqrt{1 + 
 2 \delta(\mathbf{z})}} \; \geq \; {\beta \over \sqrt{1+ 
 2 \beta}},
\end{equation*}
since at the corrector stage we have $\delta(\mathbf{z}) \geq \beta$.

Since $f(\cdot)$ is self-concordant, one damped Newton step with $\alpha = - {f'(0) \over (1+\lambda) f''(0)}$ ensures the following progress in the function value (e.g. Chapter 4 in \cite{YN-LN-2018}):
\begin{equation*}
F(\mathbf{x},\mathbf{s},\mathbf{w}) - F(\mathbf{x}_+,\mathbf{s}_+,\mathbf{w})  =  f(0) - f(\alpha) \geq
 \omega(\lambda) \; \geq \; \omega \left( {\beta \over \sqrt{1+2\beta}} \right),
\end{equation*}
where  $\mathbf{z} = (\mathbf{x},\mathbf{s},\mathbf{w})$ and $\mathbf{z}_+ = (\mathbf{x}_+,\mathbf{s}_+,\mathbf{w})$ are two consecutive iterates of the inner loop.

In view of the assumption of the lemma, we have $\Psi(\mathbf{z}) \leq \tau$. Hence, if we take $k$ steps in the corrector stage, then we have
$\tau \geq (k-1) \omega \left( {\beta \over \sqrt{1+2  \beta}} \right)$.
\endproof

\subsection{Predictor step}

Next lemma follows from a self-concordant property of a special function (see Theorem~5.1.9 in \cite{YN-LN-2018}).

\begin{lemma}{ (Lemma 5.2 in \cite{UTD_LP})}\label{technicallemma}
Consider the function $\Phi(\mathbf{b}) := - \sum\limits_{i=0}^n \ln (1+ b_i)$. If $\|\mathbf{b}\| < 1$ and $\bar{\mathbf{e}}^T \mathbf{b} = 0$, then
\begin{equation*}
\Phi(\mathbf{b}) \leq \omega_*(\|\mathbf{b}\|), \end{equation*}
where $\omega_*$ is defined in (\ref{omegas}).
\end{lemma}

In the case of LP and LCP, the definitions of the residuals $\mathbf{r}$ and $\hat{\mathbf{r}}$ and the proximity measures $\zeta_i$, $i=0,1,2$ are the same, since they are independent on the affine constraints defining the problems. Hence, the following lemma is the same as in \cite{UTD_LP}.

\begin{lemma}\label{lm-Delta}
Let $\mathbf{z} \in \mathcal{N}_{\delta}(\beta)$, with  $\beta < \frac{1}{2}$. Then, $\zeta_0(\mathbf{z}) \leq {\beta \over \sqrt{1-\beta}}$ and
\begin{align}
1-\beta  \; \leq \; & \hat r_i^2(\mathbf{z}) \; \leq \; {1 \over 1 - \beta}, \quad i=0, \dots, n, \label{eq-OBounds}\\ 
1-\beta  \; \leq \; &\frac{x_is_i}{x_i(\mathbf{w}) s_i(\mathbf{w})} \; \leq \; {1 \over 1 - \beta}, \quad i=1, \dots, n. \label{eq-XSBounds}
\end{align}
\end{lemma}
\proof{Proof.}
This proof is a slight modification of the proof of Lemma 5.3 in \cite{UTD_LP}. 
Denote by $\hat r_{\min} = \min\limits_{0 \leq i \leq n} \hat r_{i}(\mathbf{z})$ and $\hat r_{\max} = \max\limits_{0 \leq i \leq n} \hat r_i(\mathbf{z})$. 
In view of equality (\ref{res_av}), by definition of $\hat{\mathbf{r}}(\mathbf{z})$ in (\ref{t}), we have $\hat r_{\min} \leq  1 \leq \hat r_{\max}$. From the assumption $\zeta_0^2(\mathbf{z})  \leq \beta \zeta_1(\mathbf{z})$ and definition of $\delta$ in (\ref{def-delta}), we have
\begin{equation}\label{zeta_1}
\zeta_0^2(\mathbf{z})  =  \sum\limits_{i=0}^n \left( \hat r_i(\mathbf{z}) - \frac{1}{\hat r_i(\mathbf{z})} \right)^2 \; \leq \; \beta \left[ \sum\limits_{i=0}^n \left(   1 - \frac{1}{ \hat r_i(\mathbf{z})} \right)^2 \right]^{1/2} \; \leq \; \beta \frac{1}{\hat r_{\min}} \zeta_0(\mathbf{z}) .
\end{equation}
Hence, $0 \leq \frac{1}{\hat r_{\min}} - \hat r_{\min}  \leq \frac{\beta}{\hat r_{\min}} $, and we get $\hat r_i^2(\mathbf{z}) \geq 1-\beta$, for $i=0,\ldots,n$.  Thus, by (\ref{zeta_1}), we have $\zeta_0(\mathbf{z}) \leq {\beta \over \sqrt{1-\beta}}$. Therefore, $0 \leq \hat r_{\max} -\frac{1}{\hat r_{\max}} \leq \frac{\beta}{\sqrt{1-\beta}} $ and we get 
$\hat r_i^2(\mathbf{z}) \leq \frac{1}{1-\beta}$. 

For proving the inequalities (\ref{eq-XSBounds}), note that from (\ref{eq-OBounds}), we have
\begin{equation*}
x_is_i-v_i^2  =  r_i(\mathbf{z}) \; \geq \; (1-\beta) \rho(\mathbf{w}) \; = \; (1-\beta)\left[x_i(\mathbf{w})s_i(\mathbf{w}) - v_i^2\right], \quad i=1,\ldots,n.
\end{equation*}
Thus, $x_is_i \geq (1-\beta)x_i(\mathbf{w}) s_i(\mathbf{w})+ \beta v_i^2 \geq 
(1-\beta)x_i(\mathbf{w})s_i(\mathbf{w})$. The remaining inequalities can be proven in the same way.
\endproof

The following lemma is slightly different from that one in \cite{UTD_LP}, since now we fix the value of $\tau$ as a function of $\beta$. Therefore, we can have a smaller upper bound on $\beta$.

\begin{lemma}{(Lemma 5.4 in \cite{UTD_LP})}\label{neighbourhoods}
Let $\mathbf{z} \in \mathcal{N}_{\delta}(\beta)$ with $\beta < \frac{1}{3}$ and let $\tau = \omega_*\left(\frac{2\beta}{1-\beta}\right)$. Then, the predictor step length $\alpha_{p} \in (0,1)$ is well defined and $\Psi(\mathbf{z}(\alpha_{p})) = \tau$.  
\end{lemma}

\proof{Proof.}Note that the initial value $\Psi(\mathbf{z})$ can be represented as follows:
\begin{equation}
    \Psi(\mathbf{z})  =  - \sum\limits_{i=0}^n \ln\left(\hat r_i^2(\mathbf{z})\right), \; i = 0, \dots, n.
\end{equation}
Let us use Lemma \ref{technicallemma} for $\mathbf{b} = \hat{\mathbf{r}}^2(\mathbf{z}) -\mathbf{\bar{e}}$. We need to verify if the assumptions of the lemma hold. 
Since $\delta(\mathbf{z}) \leq \beta < \frac{1}{3}$, from Lemma \ref{lm-Delta}, we have
\begin{equation}\label{eq-Rho1}
\zeta_2(\mathbf{z})  \leq 
{\beta \over 1 - \beta} \; < \; 1,
\end{equation}
where $\zeta_2(\mathbf{z})$ is defined in (\ref{def-Rho}).
On the other hand, from (\ref{normt}) we have $\mathbf{\bar{e}}^T \left(\hat{ \mathbf{r}}^2(\mathbf{z}) - \bar{\mathbf{e}}\right) = 0$.
Then, using Lemma \ref{technicallemma}, we obtain
\begin{equation*}
\Psi(\mathbf{z}) =  \Phi\left(\hat{ \mathbf{r}}^2(\mathbf{z})-\mathbf{\bar{e}}\right) \; \leq \omega_*(\zeta_2(\mathbf{z})) \; \leq \; \omega_*\left({\beta \over 1 - \beta}\right) \; < \;\omega_*\left(\frac{2\beta}{1-\beta}\right) = \tau,
\end{equation*}
where we used the assumption of the lemma and monotonicity of the function $\omega^*$. Thus, $\mathbf{z}$ is in the interior of $\mathcal{N}_{\Psi}(\tau)$. Therefore, there exists a positive predictor step keeping our point in the neighbourhood $\mathcal{N}_{\Psi}(\tau)$. Using  (\ref{def-FProx}) and the barrier property of $F$, we have $\Psi(\mathbf{z}(\alpha_{p})) = \tau$.
\endproof

Let us derive an exact expression for the proximity measure $\Psi(\cdot)$ along the predictor direction. Denote
\begin{equation*}
\mathbf{d}(\alpha) = \mathbf{r}(\mathbf{z}(\alpha)) - \rho(\mathbf{w}(\alpha))\bar{\mathbf{e}} \in \mathbb{R}^{n+1}.
\end{equation*}

\begin{lemma}\label{growth_psi} 
Consider the predictor search direction $\Delta \mathbf{z} = (\Delta \mathbf{x},\Delta \mathbf{s},  -\mathbf{w} )$ defined by system (\ref{ASD}) with $\mathbf{a} \in \mathbb{R}^n$ and let 
$\mathbf{z}(\alpha) = \left(\mathbf{x}(\alpha),\mathbf{s}(\alpha),\mathbf{w}(\alpha)\right)$ ${= (\mathbf{x}+\alpha \Delta \mathbf{x},\mathbf{s}+\alpha \Delta \mathbf{s},  (1-\alpha)\mathbf{w})}$ with $\alpha \in (0,1).$  Then,
\begin{equation*}
\Psi(\mathbf{z}(\alpha)) \; = \; - \sum\limits_{i=0}^n \ln \left(1+ {1 \over \rho(\mathbf{w}(\alpha))} d_i(\alpha)\right),
\end{equation*}
and $\mathbf{d}(\alpha) = \mathbf{r}(\mathbf{z}) - \rho(\mathbf{w}) \bar{\mathbf{e}} + \alpha \mathbf{h}(\mathbf{z}) + \alpha^2 \mathbf{g}(\mathbf{z})$, with $\bar{\mathbf{e}} = (1,\ldots,1)^T \in \mathbb{R}^{n+1}$
and
\begin{align*}
h_0(\mathbf{z}) &= -v_0 - \mathbf{e}^T \mathbf{a} + \rho(\mathbf{w}) - \frac{\|\mathbf{v}\|^2}{n+1}, \quad
g_0(\mathbf{z})  =  \frac{\|\mathbf{v}\|^2}{n+1} - \Delta \mathbf{x}^T \Delta \mathbf{s},\\
h_i(\mathbf{z}) &= a_i + 2 v_i^2 + \rho(\mathbf{w}) - \frac{\|\mathbf{v}\|^2}{n+1}, \quad g_i(\mathbf{z})  =  \Delta x_i \Delta s_i  -v_i^2  + \frac{\|\mathbf{v}\|^2}{n+1}, \quad i=1, \dots, n. 
\end{align*}
\end{lemma}
\proof{Proof.}
Using (\ref{def-FProx}), we have
\begin{equation}\label{Psizalpha}
\Psi(\mathbf{z}(\alpha))  =  - \sum\limits_{i=0}^n \ln r_i(\mathbf{z}(\alpha)) + (n+1) \ln \rho(\mathbf{w}(\alpha)) = - \sum\limits_{i=0}^n \ln\left(1 + \frac{d_i(\alpha)}{\rho(\mathbf{w}(\alpha))}\right).
\end{equation}
Note that
\begin{align}\label{eq-BProg}
 \rho(\mathbf{w}(\alpha)) & =  {1 \over n+1} \Big((1-\alpha) v_0  -  \| (1-\alpha)^2 \mathbf{v}\|^2 \Big) 
 = {1-\alpha \over n+1}\left( v_0 - \| \mathbf{v} \|^2 + \alpha \| \mathbf{v} \|^2\right) \nonumber \\
 &=(1-\alpha)\left( \rho(\mathbf{w}) + {\alpha \over n+1} \| \mathbf{v} \|^2\right).
\end{align}
On the other hand, using (\ref{ASD}), for $i= 1, \ldots,n$, we get 
\begin{align}\label{rizalpha}
r_i(\mathbf{z}(\alpha)) & =  x_i(\alpha) s_i(\alpha)- (v_i(\alpha))^2 \; = \; (x_i + \alpha \Delta x_i)(s_i + \alpha \Delta s_i) - (1-\alpha)^2 v_i^2 \nonumber\\
& =  r_i(\mathbf{z}) + \alpha\left( 2v_i^2 + a_i \right) +  \alpha^2 \left(\Delta x_i \Delta s_i - v_i^2 \right).
\end{align}

\noindent Thus, combining (\ref{eq-BProg}) with (\ref{rizalpha}), we have for $i=1,\ldots,n$
\begin{align}\label{Aialpha}
d_i(\alpha) &= r_i(\mathbf{z}(\alpha)) - \rho(\mathbf{w}(\alpha)) \nonumber\\
& =   r_i(\mathbf{z}) - \rho(\mathbf{w}) + \alpha \left(a_i + 2v_i^2 + \rho(\mathbf{w}) - \frac{\|\mathbf{v}\|^2}{n+1} \right) + \alpha^2 \left(\Delta x_i \Delta s_i - v_i^2 + \frac{\|\mathbf{v}\|^2}{n+1} \right).
\end{align}
Similarly, using (\ref{ASD})
\begin{equation}\label{r0alpha}
r_0(\mathbf{z}(\alpha))  =   (1-\alpha) v_0   - \left( \mathbf{s} + \alpha \Delta \mathbf{s}\right)^T \left(\mathbf{x}+ \alpha \Delta \mathbf{x} \right)
 =   r_0(\mathbf{z}) -  \alpha \left(v_0 + \mathbf{e}^T \mathbf{a}\right)  - \alpha^2 \Delta \mathbf{x}^T \Delta \mathbf{s}.
\end{equation}
Hence, from (\ref{eq-BProg}) and (\ref{r0alpha}), we get
\begin{align}\label{A0alpha}
d_0(\alpha) &= r_0(\mathbf{z}(\alpha)) - \rho(\mathbf{w}(\alpha)) \nonumber\\
&=  r_0(\mathbf{z}) - \rho(\mathbf{w}) - \alpha \left( v_0 + \mathbf{e}^T \mathbf{a} - \rho(\mathbf{w}) + \frac{\|\mathbf{v}\|^2}{n+1} \right)+ \alpha^2 \left(\frac{\|\mathbf{v}\|^2}{n+1} -  \Delta \mathbf{x}^T \Delta \mathbf{s}\right). 
\end{align}
From (\ref{Psizalpha}), (\ref{Aialpha}) and (\ref{A0alpha}), we obtain the result of the lemma. 
 \endproof

\begin{lemma}\label{lem:d_alpha_beta} 
Let $\mathbf{z} \in \mathcal{N}_{\delta}(\beta)$ with $\beta < \frac{1}{3}$ and $\tau = \omega_*({2\beta \over 1- \beta})$. Let the predictor search direction $\Delta\mathbf{z}= (\Delta\mathbf{x},\Delta\mathbf{s},-\mathbf{w})$ be defined by system (\ref{ASD}) with $\mathbf{a}\in\mathbb{R}^n$.  
Then, for the predictor step length $\alpha_{p}$, the following inequality holds
\begin{equation*}
\frac{2\beta}{1-\beta}\rho(\mathbf{w}) (1-\alpha_{p})\le \left\|\mathbf{d}(\alpha_{p})\right\|.    
\end{equation*}
\end{lemma}
\proof{Proof.}
We consider the vectors $\mathbf{g}(\mathbf{z})$ and $\mathbf{h}(\mathbf{z})$  defined in Lemma~\ref{growth_psi}. By their definitions, we have $\mathbf{g}(\mathbf{z})^T \mathbf{\bar{e}} = 0$ and $\mathbf{h}(\mathbf{z})^T\mathbf{\bar{e}}=0$. Moreover, 
$ \left(\mathbf{r}(\mathbf{z}) - \rho(\mathbf{w}) \bar{\mathbf{e}} \right)^T \bar{\mathbf{e}}= 0$ by  (\ref{def-Rho}). 
Therefore, $\bar{\mathbf{e}}^T \mathbf{d}(\alpha)=0$. 

Assume, that 
\begin{equation}\label{assume:norm_d_alpha}
\frac{\|\mathbf{d}(\alpha_{p})\|}{\rho(\mathbf{w}(\alpha_{p}))}<\frac{2\beta}{1-\beta}.
\end{equation}
The upper bound is less than one, so by Lemma \ref{technicallemma} with $\mathbf{b} = \frac{\mathbf{d}(\alpha_{p})}{\rho(\mathbf{w}(\alpha_{p}))}$, and we get
\begin{equation*}
\Psi(\mathbf{z}(\alpha_{p}))  = \Phi\left(\frac{\mathbf{d}(\alpha_{p})}{\rho(\mathbf{w}(\alpha_{p}))}\right) \; \leq \;
\omega_*\left(\frac{\|\mathbf{d}(\alpha_{p})\|}{\rho(\mathbf{w}(\alpha_{p}))}\right).
\end{equation*}
But this contradicts to our assumption (\ref{assume:norm_d_alpha}) because of the choice of $\tau$, monotonicity of $\omega_*$, and equality $\Psi(\mathbf{z}(\alpha_{p})) = \tau$ (see Lemma \ref{neighbourhoods}).

Finally, from (\ref{eq-BProg}) we have $\rho(\mathbf{w}(\alpha)) \geq  (1-\alpha) \rho(\mathbf{w})$. Hence,
\begin{equation*}
\frac{2\beta}{1-\beta}<\frac{\|\mathbf{d}(\alpha_{p})\|}{\rho(\mathbf{w}(\alpha_{p}))}\le\frac{\|\mathbf{d}(\alpha_{p})\|}{(1-\alpha_{p})\rho(\mathbf{w})},
\end{equation*}
and we get the result.
\endproof

\begin{corollary}\label{cor:1-alpha_p}
Let $\mathbf{z} \in \mathcal{N}_{\delta}(\beta)$ with $\beta < \frac{1}{3}$ and $\tau = \omega_*\left({2\beta \over 1- \beta}\right)$. Let the predictor search direction $\Delta\mathbf{z}= (\Delta\mathbf{x},\Delta\mathbf{s},-\mathbf{w})$ be defined by system (\ref{ASD}) with $\mathbf{a}\in\mathbb{R}^n$. 
Then, for the predictor step length $\alpha_{p}$, the following inequality holds
\begin{equation*}
\frac{\beta }{1-\beta} \rho(\mathbf{w}) (1-\alpha_{p})\le 
\alpha_{p}\|\mathbf{r}(\mathbf{z}) -\rho(\mathbf{w})\bar{\mathbf{e}} + \mathbf{h}(\mathbf{z})\|+\alpha_{p}^2\left(\|\mathbf{v}\|^2+\sqrt{(\Delta\mathbf{x}^T\Delta\mathbf{s})^2+\|\Delta\mathbf{x}\Delta\mathbf{s}\|^2}\right).    
\end{equation*}
\end{corollary} 
\proof{Proof.}
From Lemma \ref{growth_psi}, we have
\begin{equation*}
\|\mathbf{d}(\alpha)\|\leq  (1-\alpha) \| \mathbf{r}(\mathbf{z}) - \rho(\mathbf{w}) \bar{\mathbf{e}} \| + \alpha \|\mathbf{r}(\mathbf{z}) - \rho(\mathbf{w}) \bar{\mathbf{e}} + \mathbf{h}(\mathbf{z}) \| +
\alpha^2 \| \mathbf{g}(\mathbf{z}) \|.   
\end{equation*}
By the definition of $\zeta_2$ in (\ref{def-Rho}) and (\ref{eq-Rho1}), we get an estimation for the first term, that is,
\begin{equation}\label{est:d_alpha_term1}
    \| \mathbf{r}(\mathbf{z}) - \rho(\mathbf{w}) \bar{\mathbf{e}} \|\le \frac{\beta}{1-\beta}\rho(\mathbf{w}).
\end{equation}
Combining it with the lower bound on the norm of $\mathbf{d}(\alpha_p)$ in Lemma \ref{lem:d_alpha_beta}, we have
\begin{equation}\label{estim:beta_g}
\frac{\beta}{1-\beta}\rho(\mathbf{w})(1-\alpha_p)\le \alpha_p {\|\mathbf{r}(\mathbf{z}) - \rho(\mathbf{w})\bar{\mathbf{e}} + \mathbf{h}(\mathbf{z}) \|} +
\alpha_p^2{\|\mathbf{g}(\mathbf{z})\|}.
\end{equation}

To give an upper bound on the last term $\| \mathbf{g}(\mathbf{z}) \|$, we consider the decomposition $\mathbf{g}(\mathbf{z}) = \mathbf{g}^{v}(\mathbf{z}) + \mathbf{g}^{\Delta}(\mathbf{z})$,  where
\begin{align*}
g_0^v(\mathbf{z}) & =   {1 \over n+1} \|\mathbf{v} \|^2, \quad
g_i^v(\mathbf{z})  = {1 \over n+1} \|\mathbf{v} \|^2  - v_i^2, \quad i=1, \ldots, n, \\ \nonumber\\
g_0^{\Delta}(\mathbf{z}) & =   -\Delta \mathbf{x}^T \Delta \mathbf{s}, \qquad g_i^{\Delta}(\mathbf{z})  =  \Delta x_i \Delta s_i, \quad i = 1, \dots, n.
\end{align*}
Since 
\begin{equation*}
\| \mathbf{g}^v(\mathbf{z}) \|^2 = {1 \over n+1} \|  \mathbf{v} \|^4 - {2 \over n+1} \| \mathbf{v} \|^4 + \sum\limits_{i=1}^n v_i^4 \; \leq \; \| \mathbf{v} \|^4,
\end{equation*}
we have
\begin{equation*}
\| \mathbf{g}(\mathbf{z}) \| \leq \| \mathbf{g}^v(\mathbf{z}) \| + \| \mathbf{g}^{\Delta}(\mathbf{z}) \| \le  \| \mathbf{v} \|^2 +  \sqrt{\left(\Delta \mathbf{x}^T \Delta \mathbf{s}\right)^2 + \|\Delta \mathbf{x} \Delta \mathbf{s}\|^2}.
\end{equation*}
This, together with (\ref{estim:beta_g}), completes the proof.
\endproof

\begin{remark}
In the decomposition of $\mathbf{d}(\alpha)$ in Lemma \ref{growth_psi}, only the term $\mathbf{h}(\mathbf{z})$ depends on $\mathbf{a}(\mathbf{z})$, not $\mathbf{g}(\mathbf{z})$. Note that $\mathbf{a}(\mathbf{z})$ appears in the right-hand side of (\ref{ASD}), i.e. it defines the particular search direction. Hence, the change of the functional proximity measure $\Psi$ during the predictor step depends on the used search direction only through the function $\mathbf{h}(\mathbf{z})$.
We have the following two special cases.
\begin{itemize}
\item In case of the UTD IPA defined by (\ref{UTD}), we have
\begin{equation}\label{hutd}
 \mathbf{h}(\mathbf{z}) = \mathbf{0}.
\end{equation}
\item In case of the AC IPA defined by (\ref{autocorr1}), we have
\begin{equation}\label{hac}
 \mathbf{h}(\mathbf{z}) =\rho(\mathbf{w})\bar{\mathbf{e}} - \mathbf{r}(\mathbf{z}).
\end{equation}
\end{itemize}
\end{remark}

Let us give a lower bound on the predictor step length $\alpha_{p}$. Following \cite{ParLCP},
we can measure the distance from the vector $\mathbf{w} \in \hbox{{\rm dom}} \, \phi$ to the boundary of $\hbox{{\rm dom}} \, \phi$ along the direction pointing to the origin by the following coefficient:
\begin{equation*}
\underline{\alpha}(\mathbf{w}) =  \frac{v_0}{\|\mathbf{v}\|^2}> 1.
\end{equation*}
This is a particular variant of the definition in \cite{ParLCP} with $\mathbf{w}^* = \mathbf{0}$ (in this paper, we deal with LCPs instead of WLCPs).

Let us give a lower bound for the predictor step $\alpha_{p}$ by estimating $\Delta \mathbf{x} \Delta \mathbf{s}$ by Lem\-ma~\ref{search-dir}. For that, we need to bound first
the norm of $\frac{\mathbf{a}}{\sqrt{\mathbf{x}\mathbf{s}}}$ using the coefficient $\underline{\alpha}(\mathbf{w})$. 

\begin{lemma}\label{gtriangle}
Let $\mathbf{z} \in \mathcal{N}_{\delta}(\beta)$ with $\beta < \frac{1}{3}$. If $\mathbf{a}_{p}(\mathbf{z})=\mathbf{a}_{ut}(\mathbf{z})$ be defined by (\ref{UTD}) or $\mathbf{a}_{p}(\mathbf{z})=\mathbf{a}_{ac}(\mathbf{z})$ be defined by (\ref{autocorr1}), then we have
\begin{equation*}
\left\| \frac{\mathbf{a}_{p}(\mathbf{z})}{\sqrt{\mathbf{x}\mathbf{s}}}\right\|^2 \leq 
\frac{\|\mathbf{v}\|^2}{1-\beta} \frac{\underline{\alpha}(\mathbf{w})^2}{\underline{\alpha}(\mathbf{w})-1}.
\end{equation*}
\end{lemma}

\proof{Proof.}
From (\ref{eq-XOpt}), we have
\begin{equation}\label{xwswtrans}
\mathbf{x}(\mathbf{w})^T\mathbf{s}(\mathbf{w}) = \frac{nv_0 + \|\mathbf{v}\|^2}{n+1}
\end{equation}
and for all $i=1,\ldots,n$
\begin{equation}\label{xwswtrans2}
x_i(\mathbf{w})s_i(\mathbf{w}) = v_i^2 + \rho(\mathbf{w}) \geq  \rho(\mathbf{w}).
\end{equation}
First, we consider $\mathbf{a}_{p}(\mathbf{z})=\mathbf{a}_{ut}(\mathbf{z})$. Then, using Lemma \ref{lm-Delta}, we have
\begin{align}\label{a_dir_ut}
\left\| \frac{\mathbf{a}_{ut}(\mathbf{z})}{\sqrt{\mathbf{x}\mathbf{s}}}\right\|^2 
& \leq {1 \over 1-\beta} \sum\limits_{i=1}^n {1 \over x_i(\mathbf{w}) s_i(\mathbf{w})} \left[ 2 x_i(\mathbf{w}) s_i(\mathbf{w})  - {1 \over n+1} v_0 \right]^2 \nonumber\\
& = {1 \over 1-\beta} \left[ 4  \mathbf{x}(\mathbf{w})^T\mathbf{s}(\mathbf{w}) - 4 {n \over n+1} v_0 + {(v_0)^2 \over (n+1)^2}  \sum\limits_{i=1}^n {1 \over x_i(\mathbf{w}) s_i(\mathbf{w}) } \right].
\end{align}

For auto-corrector direction  $\mathbf{a}_{p}(\mathbf{z})=\mathbf{a}_{ac}(\mathbf{z})$, by Lemma \ref{lm-Delta}, we get the same bound: 
\begin{align}\label{a_dir_ac}
\left\| \frac{\mathbf{a}_{ac}(\mathbf{z})}{\sqrt{\mathbf{x}\mathbf{s}}}\right\|^2 
&\leq \sum\limits_{i=1}^n \frac{1}{x_is_i}
\left[x_i(\mathbf{w}) s_i(\mathbf{s}) - \frac{v_0}{n+1} \right]^2 
+\sum_{i=1}^n x_i s_i + 2 \sum_{i=1}^n \left[x_i(\mathbf{w}) s_i(\mathbf{w}) - \frac{v_0}{n+1}\right]  \nonumber\\
&\leq {1 \over 1-\beta} \sum\limits_{i=1}^n {1 \over x_i(\mathbf{w}) s_i(\mathbf{w})} \left[  x_i(\mathbf{w}) s_i(\mathbf{w})  - {v_0 \over n+1} \right]^2 +  \frac{3-2\beta}{1-\beta}\sum_{i=1}^n  
 x_i(\mathbf{w})s_i(\mathbf{w})  -\frac{2nv_0}{n+1}\nonumber\\
&= {1 \over 1 - \beta} \left[ (4-2\beta)  \mathbf{x}(\mathbf{w})^T\mathbf{s}(\mathbf{w}) - 2  {n \over n+1} v_0 + {(v_0)^2 \over (n+1)^2}  \sum\limits_{i=1}^n {1 \over x_i(\mathbf{w}) s_i(\mathbf{w}) } \right]   -\frac{2nv_0}{n+1}\nonumber\\
&\leq {1 \over 1-\beta} \left[ 4  \mathbf{x}(\mathbf{w})^T\mathbf{s}(\mathbf{w}) - 4 {n \over n+1} v_0 + {(v_0)^2 \over (n+1)^2}  \sum\limits_{i=1}^n {1 \over x_i(\mathbf{w}) s_i(\mathbf{w}) } \right],
\end{align}
where in the last inequality we used $\mathbf{x}(\mathbf{w})^T\mathbf{s}(\mathbf{w}) \geq \frac{nv_0}{n+1}$ obtained from (\ref{xwswtrans}).
We further estimate the bounds obtained in (\ref{a_dir_ut}) and (\ref{a_dir_ac}) by using (\ref{xwswtrans}) and (\ref{xwswtrans2}) as follows
\begin{equation*}
\left\| \frac{\mathbf{a}_{p}(\mathbf{z})}{\sqrt{\mathbf{x}\mathbf{s}}}\right\|^2 
\leq {1 \over 1-\beta} \left[ {4 \| \mathbf{v} \|^2 \over n+1} + {n v_0^2 \over (n+1)^2 \rho(\mathbf{w})} \right]
=  {\| \mathbf{v} \|^2 \over (1 - \beta)(n+1)} \left[ 4 + {n  \underline{\alpha}(\mathbf{\mathbf{w}})^2 \over \underline{\alpha}(\mathbf{\mathbf{w}}) - 1} \right].
\end{equation*}
Using that $\underline{\alpha}(\mathbf{w})>1$, we have
\begin{equation}\label{alpha_under_4}
4 \leq \frac{\underline{\alpha}(\mathbf{\mathbf{w}})^2}{ \underline{\alpha}(\mathbf{\mathbf{w}}) - 1},
\end{equation}
which completes the proof.
\endproof

Let us show that for both our predictor directions, we have the same lower bounds for the predictor step length $\alpha_{p}$.

\begin{lemma}\label{lm-HMin}
Let $\mathbf{z} \in \mathcal{N}_{\delta}(\beta)$ with $\beta < \frac{1}{3}$ and $\tau = \omega_*\left({2\beta \over 1- \beta}\right)$. If $\mathbf{a}_{p}(\mathbf{z})=\mathbf{a}_{ut}(\mathbf{z})$ be defined by (\ref{UTD}) or $\mathbf{a}_{p}(\mathbf{z})=\mathbf{a}_{ac}(\mathbf{z})$ be defined by (\ref{autocorr1}), then the predictor step length $\alpha_{p}$ satisfies the following condition:
\begin{equation}\label{eq-HMin}
 0 \; < \;  {\beta \over 1 - \beta}  \leq  {2\beta \over 1 - \beta}  \alpha_{p} + 
   (n+1) \left( {\underline{\alpha}(\mathbf{w}) \over \underline{\alpha}(\mathbf{w}) - 1} \right)^2 \alpha_{p}^2.
 \end{equation}
\end{lemma}
\proof{Proof.}
From Lemma \ref{search-dir}
and Lemma \ref{gtriangle}, we derive
\begin{equation*}
 \sqrt{\left(\Delta \mathbf{x}^T \Delta \mathbf{s}\right)^2 + \|\Delta \mathbf{x} \Delta \mathbf{s}\|^2} \leq \frac{1}{2} \left\| \frac{\mathbf{a}_{p}(\mathbf{z})}{\sqrt{\mathbf{x}\mathbf{s}}} \right\|^2 \leq \frac{1}{2}\frac{\|\mathbf{v}\|^2}{1-\beta} \frac{\underline{\alpha}(\mathbf{w})^2}{\underline{\alpha}(\mathbf{w})-1}.
 \end{equation*}
Therefore by (\ref{alpha_under_4}), the identity $\| \mathbf{v} \|^2 = {n+1 \over \underline{\alpha}(\mathbf{w}) - 1}\rho(\mathbf{w})$, and since $\beta < {1 \over 3}$, we have
\begin{equation}\label{normgz} 
\| \mathbf{v} \|^2 \Big[ 1 +\frac{1}{2} { \underline{\alpha}(\mathbf{w})^2 \over (1 - \beta)(\underline{\alpha}(\mathbf{w}) - 1)} \Big] \; <  (n+1) \left( {\underline{\alpha}(\mathbf{w}) \over \underline{\alpha}(\mathbf{w}) - 1} \right)^2 \rho(\mathbf{w}).
\end{equation}

Combining (\ref{estim:beta_g}) and (\ref{normgz}), 
\begin{equation}\label{estim:pred_alpa_from_d}
 (1-\alpha_p)\frac{\beta}{1-\beta} \le  {\alpha_p} \frac{\|\mathbf{r}(\mathbf{z}) - \rho(\mathbf{w})\bar{\mathbf{e}} + \mathbf{h}(\mathbf{z}) \|}{\rho(\mathbf{w})} +
\alpha_p^2 (n+1) \left( {\underline{\alpha}(\mathbf{w}) \over \underline{\alpha}(\mathbf{w}) - 1} \right)^2.
\end{equation}

In the case of $\mathbf{a}_{p}(\mathbf{z}) = \mathbf{a}_{ut}(\mathbf{z})$, the vector $\mathbf{h}$ is given by (\ref{hutd}), so 
$\mathbf{r}(\mathbf{z}) - \rho(\mathbf{w})\bar{\mathbf{e}} + \mathbf{h}(\mathbf{z}) = \mathbf{r}(\mathbf{z}) - \rho(\mathbf{w})\bar{\mathbf{e}}$.
While in the case of $\mathbf{a}_{p}(\mathbf{z}) = \mathbf{a}_{ac}(\mathbf{z})$, the vector $\mathbf{h}$ is given by (\ref{hac}), and so  $\mathbf{r}(\mathbf{z}) - \rho(\mathbf{w})\bar{\mathbf{e}} + \mathbf{h}(\mathbf{z}) = \mathbf{0}$.
Therefore, in both cases $\|\mathbf{r}(\mathbf{z}) - \rho(\mathbf{w})\bar{\mathbf{e}} + \mathbf{h}(\mathbf{z})\|\le \frac{\beta}{1-\beta}\rho(\mathbf{w})$ by (\ref{est:d_alpha_term1}).
Substituting this upper bound in (\ref{estim:pred_alpa_from_d}), we complete the proof.
\endproof

The quadratic inequality (\ref{eq-HMin}) for the predictor step length ensures that we can give a lower bound on $\alpha_p$ of magnitude $\frac{1}{\sqrt{n}}$. 

\begin{corollary}\label{lm-HMin2}
Let $\mathbf{z} \in \mathcal{N}_{\delta}(\beta)$ with $\beta < \frac{1}{3}$ and $\tau = \omega_*\left({2\beta \over 1 - \beta} \right)$. If $\mathbf{a}_{p}(\mathbf{z})=\mathbf{a}_{ut}(\mathbf{z})$ be defined by (\ref{UTD}) or $\mathbf{a}_{p}(\mathbf{z})=\mathbf{a}_{ac}(\mathbf{z})$ be defined (\ref{autocorr1}), then we have $$
\alpha_{p}  > \frac{1}{2}\sqrt{\frac{\beta}{1-\beta}} \frac{\underline{\alpha}(\mathbf{w})-1}{\underline{\alpha}(\mathbf{w})} \frac{1}{\sqrt{n}}.$$
\end{corollary}
\proof{Proof.}
By Lemma \ref{lm-HMin}, we know
\[
1\le 2\alpha_p + t\alpha_p^2,
\]
where $t=(n+1)\left( {\underline{\alpha}(\mathbf{w}) \over \underline{\alpha}(\mathbf{w}) - 1}\right)^2\frac{1-\beta}{\beta}$. Since $n\ge2$, we have $3\le n+1\le\frac{3}{2} n$. Furthermore, $\frac{1-\beta}{\beta}>2$, thus
\[6<t\le\frac{3}{2}n\left( {\underline{\alpha}(\mathbf{w}) \over \underline{\alpha}(\mathbf{w}) - 1}\right)^2\frac{1-\beta}{\beta}.\]
On the other hand, the above quadratic inequality on $\alpha_p$ means that
\begin{align*}
    \alpha_p&\ge \frac{-1+\sqrt{1+t}}{t}=\frac{1}{\sqrt{t}}\left(\sqrt{\frac{1}{t}+1}-\sqrt{\frac{1}{t}}\right)\\
    &> \sqrt{\frac{2}{3}}\frac{1}{\sqrt{n}} {\underline{\alpha}(\mathbf{w})-1 \over \underline{\alpha}(\mathbf{w})}\sqrt{\frac{\beta}{\beta-1}}\left(\sqrt{\frac{1}{6}+1}-\sqrt{\frac{1}{6}}\right)>
    \frac{1}{2\sqrt{n}} {\underline{\alpha}(\mathbf{w})-1 \over \underline{\alpha}(\mathbf{w})}\sqrt{\frac{\beta}{\beta-1}},
\end{align*}
where we used the monotonicity function $f(x)=\sqrt{x+1}-\sqrt{x}$.
\endproof

As in \cite{ParLCP}, we can establish a global convergence rate of the sequence generated by Algorithm \ref{alg:PTS} in terms of the following merit function 
\begin{equation*}
 \mu^*(\mathbf{w}) = \frac{\underline{\alpha}(\mathbf{w})}{\underline{\alpha}(\mathbf{w})-1}v_0 = \frac{v_0^2}{v_0-\|\mathbf{v}\|^2} \geq \; v_0 \; \geq \; \| \mathbf{v} \|^2.
\end{equation*}

Let us give an upper bound for the new value of the merit function after the predictor step.

\begin{lemma}{ (Lemma 5.7 in \cite{UTD_LP})}\label{lem:important1}
Let $\alpha_{p} \in (0, 1)$ be a feasible predictor step of PTS IPA.
If $\alpha_{p} \geq \gamma \, \frac{\underline{\alpha}(\mathbf{w}) - 1}{\underline{\alpha}(\mathbf{w})}$ with some $\gamma \in (0,1)$, then $\mu^*(\mathbf{w}(\alpha_{p})) < \frac{1}{1+\gamma} \, \mu^*(\mathbf{w})$.
\end{lemma}
\proof{Proof.}
For the sake of simplicity, we use the notation $\alpha=\alpha_{p}$. By the assumption of the lemma and using that $\underline{\alpha}(\mathbf{w})> 1$, we have
\begin{equation*}
\frac{1}{1-\alpha}=\frac{1+\alpha}{1-\alpha^2}> 1+\alpha\ge 1+\gamma \,\frac{\underline{\alpha}(\mathbf{w})-1}{\underline{\alpha}(\mathbf{w})}.
\end{equation*}

Since $\underline{\alpha}(\mathbf{w}(\alpha))=\frac{1}{1-\alpha}\,\underline{\alpha}(\mathbf{w})$,
\begin{align*}
    \frac{\mu^*(\mathbf{w}(\alpha))}{\mu^*(\mathbf{w})}&= \frac{\frac{1}{1-\alpha}\,\underline{\alpha}(\mathbf{w})}{\frac{1}{1-\alpha}\,\underline{\alpha}(\mathbf{w})-1}\,(1-\alpha)\,v_0\, \frac{\underline{\alpha}(\mathbf{w})-1}{\underline{\alpha}(\mathbf{w})v_0}=\frac{\underline{\alpha}(\mathbf{w})-1}{\frac{1}{1-\alpha}\,\underline{\alpha}(\mathbf{w})-1}\\
    &< \frac{\underline{\alpha}(\mathbf{w})-1}{\left(1+\gamma \,\frac{\underline{\alpha}(\mathbf{w})-1}{\underline{\alpha}(\mathbf{w})}\right)\underline{\alpha}(\mathbf{w})-1}=\frac{1}{1+\gamma},
\end{align*}
hence the lemma is proved.
\endproof

Now we are ready to give the complexity results for PTS IPA.
\begin{theorem}\label{nrouteriter}
Let $ \beta = \frac{1}{4}$, $\tau = \omega_*\left({2\beta \over 1- \beta}\right)$, and $\mathbf{z}^{(0)} = (\mathbf{x}^{(0)},\mathbf{s}^{(0)},\mathbf{w}^{(0)}) \in \mathcal{F}_z^+$  be the starting point. If $\mathbf{a}_{p}(\mathbf{z})=\mathbf{a}_{ut}(\mathbf{z})$ be defined by (\ref{UTD}) or $\mathbf{a}_{p}(\mathbf{z})=\mathbf{a}_{ac}(\mathbf{z})$ be defined by (\ref{autocorr1}), then the PTS IPA gives a feasible solution to (\ref{LCP}) with $\mathbf{x}^T\mathbf{s} \leq \varepsilon$ after at most 
\begin{equation*}
\left \lceil 4\sqrt{n} \ln {\mu^*(\mathbf{w}^{(0)}) \over \varepsilon} \right \rceil \quad \mbox{iterations.}
\end{equation*}
\end{theorem}
    
\proof{Proof.} By Corollary \ref{lm-HMin2},  $\gamma = \frac{1}{2}\sqrt{\frac{\beta}{1-\beta}} \frac{1}{\sqrt{n}} = \frac{1}{\sqrt{12}} \frac{1}{\sqrt{n}} > \frac{1}{4 \sqrt{n}}$ satisfies the condition of  Lemma \ref{lem:important1}.
Hence, 
$\mu^*(\mathbf{w}(\alpha_{p})) < {4\sqrt{n} \over 1+ 4\sqrt{n} } \mu^*(\mathbf{w})$, and after $k$ outer iterations, we have
$v_0 \leq \mu^*(\mathbf{w}) < \left(\frac{4\sqrt{n}}{1+4\sqrt{n}}\right)^k \mu^*\left(\mathbf{w}^{(0)}\right)$. 

Note that the stopping criterion of the algorithm $v^0 \leq \varepsilon$ is satisfied if
$ \left(\frac{4\sqrt{n}}{1+4\sqrt{n}}\right)^k  \mu^*(\mathbf{w}^{(0)}) \leq \varepsilon$. 
Since $\ln (1+t)  \leq t$ for $t\geq 0$ and since the lengths of corrector stages are uniformly bounded by inequality (\ref{eq-NUp}), we get the result.
 \endproof

\begin{remark}
As we have described at the beginning of Section 3, for any $(\mathbf{x}^{(0)},\mathbf{s}^{(0)}) \in 
\mathcal{F}^+$, we can always compute the corresponding $\mathbf{w}^{(0)}=\mathbf{w}(\mathbf{x}^{(0)},\mathbf{s}^{(0)}) \in \mathcal{F}_w^+$ such that $\delta\left(\mathbf{x}^{(0)},\mathbf{s}^{(0)},\mathbf{w}^{(0)}\right) = 0$. Hence, from (\ref{def-W(x)})
we have 
\[
\mu^*\left(\mathbf{w}^{(0)}\right)= \frac{\left(\left(\mathbf{x}^{(0)}\right)^T\!\mathbf{s}^{(0)}+\xi\left(\mathbf{x}^{(0)},\mathbf{s}^{(0)}\right)\right)^2}{2 \xi\left(\mathbf{x}^{(0)},\mathbf{s}^{(0)}\right)} = \mathcal{O}\left(
\left(\frac{\left(\mathbf{x}^{(0)}\right)^T\!\mathbf{s}^{(0)}}{\sqrt{\min\limits_{1 \leq i \leq n} x^{(0)}_i s^{(0)}_i}}\right)^2\right).
\]
Thus, the complexity of the algorithm is
$$ \mathcal{O}\left(\sqrt{n} \ln \left(\frac{\left(\mathbf{x}^{(0)}\right)^T\!\mathbf{s}^{(0)}}{\varepsilon \cdot \sqrt{\min\limits_{1 \leq i \leq n} x^{(0)}_i s^{(0)}_i}}\right) \right). $$
Hence, the condition measure of the starting point $\left(\mathbf{x}^{(0)},\mathbf{s}^{(0)}\right)$ appears only inside the logarithm, contrary to the older versions of target-following IPAs (e.g.
 \cite{jrt}). 
\end{remark}

\section{Analysis of local convergence}

In this section, we analyze the properties of the auto-corrector direction under the following non-degeneracy assumption.

\textbf{Assumption.}
The (\ref{LCP}) has a unique and strictly complementary solution $(\mathbf{x}^*,\mathbf{s}^*)$. We assume that $\mathbf{x}^*$ has $m$ positive components. 

Without loss of generality, these are the first $m$ components of the vector, i.e.
$$ \mathbf{x}^* = (\mathbf{x}_B^*, \mathbf{x}_N^{*}), \quad \mathbf{x}_B^* \in \mathbb{R}_+^m, \quad \mathbf{x}_N^* = \mathbf{0} \in \mathbb{R}^{n-m}.$$
Otherwise, we permute the columns and rows of the matrix $M$.
We have
$$ \mathbf{s}^* = (\mathbf{s}_B^*, \mathbf{s}_N^{*}), \quad \mathbf{s}_B^* =\mathbf{0}  \in \mathbb{R}^m, \quad \mathbf{s}_N^*  \in \mathbb{R}_{+}^{n-m}.$$

Note that from the uniqueness of the solution it follows that $M_{BB}$ is invertible.

We need to give an upper bound for the norm of $\Delta\mathbf{x}\Delta\mathbf{s}$. Lemma \ref{search-dir} is not suitable for that, since we cannot give a proper lower bound on $x_is_i$ (note that $x_i$, $i\in N$ and $s_i$, $i\in B$ tend to zero through the progress of the algorithm). Therefore, we consider another rescaling of the system (\ref{ASD}), where we split equations in accordance to positivity of coordinates of $\mathbf{x}$ and $\mathbf{s}$ in the solution. Then we can give proper lower bounds for them.

In view of the partition $(B,N)$ of the index set, we can consider the system defining the search directions (\ref{ASD}) in the following form:
\begin{align*}
M_{BB} \Delta \mathbf{x}_B + M_{BN} \Delta \mathbf{x}_N &= \Delta \mathbf{s}_B, \nonumber\\
M_{NB} \Delta \mathbf{x}_B + M_{NN} \Delta \mathbf{x}_N &= \Delta \mathbf{s}_N, \\
S_B \Delta \mathbf{x}_B + X_B \Delta \mathbf{s}_B &= \mathbf{a}_B, \\
 X_N \Delta \mathbf{s}_N + S_N \Delta \mathbf{x}_N &= \mathbf{a}_N.
\end{align*}

From the first two equations we express $\Delta \mathbf{x}_B$ and $\Delta \mathbf{s}_N$, while from the last two equations we express $\Delta \mathbf{s}_B$ and $\Delta \mathbf{x}_N$. For that, let us introduce
\begin{equation*}
\hat{M} = \left[ 
\begin{array}{cc} 
M_{BB}^{-1} \quad & \quad-M_{BB}^{-1}M_{BN} \\
M_{NB}M_{BB}^{-1} \quad &\quad M_{NN}-M_{NB}M_{BB}^{-1}M_{BN}
\end{array}
\right],
\end{equation*}
and
\begin{equation*}
D(\mathbf{x},\mathbf{s}) = \left[ \begin{array}{cc}
X_B^{-1} S_B  \quad & \quad O\\
O \quad & \quad S_N^{-1} X_N 
\end{array}\right], \quad \tilde{\mathbf{a}}(\mathbf{x},\mathbf{s}) = \left( \begin{array}{cc} X_B^{-1} \mathbf{a}_B \\
 S_N^{-1} \mathbf{a}_N \end{array} \right),
\end{equation*}
where $O$ is the null matrix of the corresponding size.
In this way, we have
\begin{equation}\label{DeltaxBDeltasN}
\left(\begin{array}{c}
\Delta \mathbf{x}_B \\
\Delta \mathbf{s}_N
\end{array}\right)
= \hat{M} 
\left(
\begin{array}{c} 
\Delta \mathbf{s}_B\\
\Delta \mathbf{x}_N
\end{array}
\right)
\end{equation}
and
\begin{equation}\label{DeltasbDeltaxn}
\left(
\begin{array}{c} 
\Delta \mathbf{s}_B\\
\Delta \mathbf{x}_N
\end{array}
\right) = 
\begin{array}{cc}
-D(\mathbf{x},\mathbf{s}) \left(\begin{array}{c}
\Delta \mathbf{x}_B \\
\Delta \mathbf{s}_N
\end{array}\right) + \tilde{{\mathbf{a}}}(\mathbf{x},\mathbf{s}),
\end{array}
\end{equation}

We introduce the following notations:
$$x_{\min}^* = \min_{1\leq i \leq m} x_i^*, \quad s_{\min}^* =\min_{m+1 \leq i \leq n} s_i^*,  $$

$$ \sigma = \min \{ x^*_{\min}, s^*_{\min}\},
\quad \kappa = \|\hat{M}\|_{\infty}, \quad \tilde{\nu}_{D}(\mathbf{x},\mathbf{s}) = \|D(\mathbf{x},\mathbf{s})\|, \quad \tilde{\nu}_{\mathbf{a}}(\mathbf{x},\mathbf{s}) = \|\tilde{\mathbf{a}}(\mathbf{x},\mathbf{s})\|.$$

Thus, $\tilde{\nu}_{D}(\mathbf{x},\mathbf{s})$ and $\tilde{\nu}_{\mathbf{a}}(\mathbf{x},\mathbf{s})$ are dependent on the actual iterate, while $\kappa$ and $\sigma$ do not depend on the actual iterate. Note that $\sigma$ is the minimal nonzero component of the unique solution $(\mathbf{x}^*, \mathbf{s}^*)$, and $$\tilde{\nu}_{D}(\mathbf{x},\mathbf{s})=\max_{1\leq i \leq n} D_{ii}(\mathbf{x},\mathbf{s}) =
\max\left\{ \frac{s_i}{x_i}, \frac{x_j}{s_j} :\; i\in B,\; j\in N \right\}. $$

From now on, we denote $D(\mathbf{x},\mathbf{s})$ by $D$, $\tilde{\mathbf{a}}(\mathbf{x},\mathbf{s})$ by $\tilde{\mathbf{a}}$, $\tilde{\nu}_{D}(\mathbf{x},\mathbf{s})$ by $\tilde{\nu}_D$ and   $\tilde{\nu}_{\mathbf{a}}(\mathbf{x},\mathbf{s})$ by $\tilde{\nu}_{\mathbf{a}}$. In what follows, we will give an estimate on $\|\Delta \mathbf{x} \Delta \mathbf{s}\|$.

\begin{lemma}\label{technical1}
Let $(\mathbf{x},\mathbf{s}) \in \mathcal{F}$. Then, we have 
\begin{equation}\label{dualgap}
\left( \begin{array}{cc}
\mathbf{x}_B^* \\
\mathbf{x}_N 
\end{array}
\right)^T \left( \begin{array}{cc}
\mathbf{s}_B \\
\mathbf{s}_N^* 
\end{array}
\right) \leq \mathbf{x}^T\mathbf{s}
\end{equation}
and
\begin{equation}\label{xb_min_xbstar}
\left\|
\left( \begin{array}{cc}
\mathbf{x}_B - \mathbf{x}_B^* \\
\mathbf{s}_N - \mathbf{s}_N^* 
\end{array}
\right)
\right\|_{\infty}
\leq \kappa \left\| \left( \begin{array}{cc} 
\mathbf{s}_B \\ \mathbf{x}_N
\end{array}
\right)
\right\|_1 \leq \frac{\kappa}{\sigma} \mathbf{x}^T\mathbf{s}.
\end{equation}
\end{lemma}
\proof{Proof.}
We have
$$ 0 \leq (\mathbf{s} - \mathbf{s}^*)^T(\mathbf{x} - \mathbf{x}^*) = \mathbf{s}^T\mathbf{x} - (\mathbf{s}^*)^T \mathbf{x} - \mathbf{s}^T\mathbf{x}^*,$$
and we obtain (\ref{dualgap}).
Using (\ref{DeltaxBDeltasN}) and $M(\mathbf{x}-\mathbf{x}^*) = \mathbf{s}-\mathbf{s}^*$,   we have
\begin{equation*}
\left( \begin{array}{cc} \mathbf{x}_B- \mathbf{x}_B^* \\ \mathbf{s}_N - \mathbf{s}_N^* \end{array}\right)= \hat{M} \left(\begin{array}{c} 
\mathbf{s}_B \\ \mathbf{x}_N\end{array}\right). 
\end{equation*}
From this we get (\ref{xb_min_xbstar}).
 \endproof

\begin{corollary}\label{technical2}
Let $(\mathbf{x},\mathbf{s}) \in \mathcal{F}$. Then, we have
\begin{equation}\label{xisi}
\left( \begin{array}{c} 
\mathbf{x}_B \\ \mathbf{s}_N\end{array}
\right) \geq \left(\sigma - \frac{\kappa \mathbf{x}^T\mathbf{s}}{\sigma}\right) \mathbf{e} \in \mathbb{R}^{2n}.
\end{equation}
\end{corollary}
\proof{Proof.}
Using  the definition of $\sigma$ and (\ref{xb_min_xbstar}), for any index $i$ we have
\begin{equation*}
\left( \begin{array}{c} 
\mathbf{x}_B \\ \mathbf{s}_N\end{array}
\right)_i \geq 
\min_i \left( \begin{array}{c} \mathbf{x}_B^* \\
\mathbf{s}_N^* \end{array} \right)_i - \left\| \left(\begin{array}{c}
\mathbf{x}_B \\ \mathbf{s}_N \end{array} \right) - \left( \begin{array}{c} \mathbf{x}_B^* \\ \mathbf{s}_N^* \end{array} \right) \right\|_{\infty} \geq 
\sigma - \frac{\kappa \mathbf{x}^T\mathbf{s}}{\sigma},
\end{equation*}
which proves the corollary.
 \endproof

 \begin{lemma}\label{rhonew}
Suppose that $\tilde{\nu}_D  \leq \frac{1}{\kappa}$. Then, we have
\begin{equation}\label{rhodelta}
\left\|  \left( \begin{array}{c} \Delta \mathbf{s}_B \\ \Delta \mathbf{x}_N \end{array}   \right)\right\| \leq \frac{\tilde{\nu}_{\mathbf{a}}}{1- \kappa \tilde{\nu}_D } \qquad\mbox{ and} \qquad \left\|  \left( \begin{array}{c} \Delta \mathbf{x}_B \\ \Delta \mathbf{s}_N \end{array}   \right)\right\| \leq \frac{\kappa \tilde{\nu_{\mathbf{a}}}}{1- \kappa \tilde{\nu}_D }.
\end{equation}  
\end{lemma}
\proof{Proof.}
From (\ref{DeltaxBDeltasN}) and (\ref{DeltasbDeltaxn}), we get
\begin{equation*}
(I+D\hat{M}) \left(
\begin{array}{c}
\Delta \mathbf{s}_B \\
\Delta \mathbf{x}_N
\end{array}
\right) = \tilde{\mathbf{a}} \qquad \text{ and }\qquad  (I+\hat{M}D) \left(
\begin{array}{c}
\Delta \mathbf{x}_B \\
\Delta \mathbf{s}_N
\end{array}
\right) = \hat{M}\tilde{\mathbf{a}},
\end{equation*}
which yields the inequalities (\ref{rhodelta}).
 \endproof

\begin{lemma}
Suppose $\mathbf{x}^T\mathbf{s} < \frac{\sigma^2}{\kappa}$. Then, we have
\begin{equation}\label{rho}
\tilde{\nu}_{D}  \leq \frac{\mathbf{x}^T\mathbf{s}}{\sigma^2 - \kappa \mathbf{x}^T\mathbf{s}} \qquad\mbox{ and} \qquad \tilde{\nu}_{\mathbf{a}}  \leq \frac{\sigma}{\sigma^2 - \kappa {\mathbf{x}^T\mathbf{s}}} \|\mathbf{a}\|.
\end{equation}
\end{lemma}
\proof{Proof.}
Using (\ref{xb_min_xbstar}) and (\ref{xisi}), we have
\begin{equation*}
\max_{i\in B} \frac{s_i}{x_i} \leq \frac{1}{\sigma-\frac{\kappa \mathbf{x}^T\mathbf{s}}{\sigma}} \frac{\mathbf{x}^T\mathbf{s}}{\sigma} = \frac{\mathbf{x}^T\mathbf{s}}{\sigma^2 - \kappa \mathbf{x}^T\mathbf{s}}
\end{equation*}
and
\begin{equation*}
\max_{j\in N}\frac{x_j}{s_j} \leq \frac{1}{\sigma-\frac{\kappa \mathbf{x}^T\mathbf{s}}{\sigma}} \frac{\mathbf{x}^T\mathbf{s}}{\sigma} = \frac{\mathbf{x}^T\mathbf{s}}{\sigma^2 - \kappa \mathbf{x}^T\mathbf{s}}.
\end{equation*}
Moreover, using (\ref{xisi}), we have
$$ \tilde{\nu}_{\mathbf{a}} =  \left\|
\left(\begin{array}{c}
X_B^{-1} \mathbf{a}_B \\
S_N^{-1} \mathbf{a}_N
\end{array} \right) \right\|\leq \frac{\sigma}{\sigma^2 - \kappa \mathbf{x}^T\mathbf{s}} \| \mathbf{a}\|,$$
which proves the lemma.
\endproof

Now we are ready to give a new upper bound for the norms of $\Delta\mathbf{x}\Delta\mathbf{s}$ as a consequence of the lower bounds for $\mathbf{x}_B$ and $\mathbf{s}_N$. 

\begin{corollary}\label{DeltaxDeltas_estimate}
Suppose that $\tilde{\nu}_D  \leq \frac{1}{\kappa}$ and $\mathbf{x}^T \mathbf{s} \leq \frac{\sigma^2}{4\kappa}$. Then, we have
\begin{equation}\label{DeltaxDeltas}
\| \Delta \mathbf{x} \Delta \mathbf{s}\|^2 \leq \frac{8 \kappa^2}{\sigma^4} \|\mathbf{a}\|^4 \quad \text{ and } \quad 
\left(\Delta \mathbf{x}^T \Delta \mathbf{s}\right)^2 \leq  \frac{8 \kappa^2}{\sigma^4} \|\mathbf{a}\|^4.
\end{equation}
\end{corollary}
\proof{Proof.}
From (\ref{rhodelta}), (\ref{rho}) and using $\mathbf{x}^T \mathbf{s} \leq \frac{\sigma^2}{4\kappa}$, we have
\begin{align*}
    \| \Delta \mathbf{x} \Delta \mathbf{s} \|^2 &=
    \left\| \left(  \begin{array}{c}
         \Delta \mathbf{x}_B \\
         \Delta \mathbf{s}_N 
    \end{array}  \right)  \left(  \begin{array}{c}
         \Delta \mathbf{s}_B \\
         \Delta \mathbf{x}_N 
    \end{array}  \right)\right\|^2 \leq 
     \left\| \left(  \begin{array}{c}
         \Delta \mathbf{x}_B \\
         \Delta \mathbf{s}_N 
    \end{array}  \right)  \right\|^2 \left\|\left(  \begin{array}{c}
         \Delta \mathbf{s}_B \\
         \Delta \mathbf{x}_N 
    \end{array}  \right)\right\|^2 \nonumber\\
    &\leq \frac{\kappa^2 \tilde{\nu}_{\mathbf{a}}^4}{(1-\kappa \tilde{\nu}_D)^2} \leq \frac{\kappa^2 \| \mathbf{a}\|^4}{\left(\sigma - \kappa \frac{\mathbf{x}^T \mathbf{s}}{\sigma} \right)^4 \left(1-\kappa \frac{\mathbf{x}^T\mathbf{s}}{\sigma^2 - \kappa \mathbf{x}^T\mathbf{s}} \right)^2} \leq \frac{8 \kappa^2}{\sigma^4} \|\mathbf{a}\|^4
\end{align*}
Thus,
\begin{equation*}
\left(\Delta \mathbf{x}^T \Delta \mathbf{s}\right)^2 \leq \left\| \left(  \begin{array}{c}
         \Delta \mathbf{x}_B \\
         \Delta \mathbf{s}_N 
    \end{array}  \right)  \right\|^2 \left\|\left(  \begin{array}{c}
         \Delta \mathbf{s}_B \\
         \Delta \mathbf{x}_N 
    \end{array}  \right)\right\|^2  \leq \frac{8 \kappa^2}{\sigma^4} \|\mathbf{a}\|^4,
    \end{equation*}
which gives the proof.
 \endproof

\subsection{Local quadratic convergence of the auto-corrector IPA}

In what follows, we analyse the local quadratic convergence of the auto-corrector IPA. In order to be able to apply Corollary \ref{DeltaxDeltas_estimate}, we assume that $v_0\le\frac{\sigma^2}{4\kappa}$. For the sake of simplicity, we will use the notation $\alpha=\alpha_p$ for the predictor step length.
In this case, $\mathbf{h}(\mathbf{z})=\rho(\mathbf{w})\bar{\mathbf{e}}-\mathbf{r}(\mathbf{z})$. Therefore, Corollary \ref{cor:1-alpha_p} and inequalities (\ref{DeltaxDeltas}) give the following bound on the predictor step length
\begin{equation}\label{estim:alpha_p-a}
\frac{\beta}{1-\beta}\rho(\mathbf{w})(1-\alpha) \leq \|\mathbf{v}\|^2 + \sqrt{\left(\Delta \mathbf{x}^T \Delta \mathbf{s}\right)^2 + \| \Delta \mathbf{x} \Delta \mathbf{s}\|^2} \leq \|\mathbf{v}\|^2 + \frac{4 \kappa}{\sigma^2} \|\mathbf{a}\|^2, 
\end{equation}
where now
\begin{equation*}
\mathbf{a} = \mathbf{a}_{ac}(\mathbf{z}) = \mathbf{a}_{ut}(\mathbf{z}) + \mathbf{a}_{c}(\mathbf{z}) = \left(\frac{\|\mathbf{v}\|^2}{n+1}  - \rho(\mathbf{w})\right) \mathbf{e} - 2\mathbf{v}^2 + \rho(\mathbf{w}) \mathbf{e} - \mathbf{r}(\mathbf{z}).
\end{equation*}
Using (\ref{est:d_alpha_term1}) and the definition of $\rho(\mathbf{w})$, we have
\begin{equation}\label{estim:cr}
\|\mathbf{a}_{c}(\mathbf{z})\|^2=\|\rho(\mathbf{w}) \mathbf{e} - \mathbf{r}(\mathbf{z})\|^2\le \frac{\beta^2}{(1-\beta)^2}\rho^2(\mathbf{w})\le
\frac{\beta^2\rho(\mathbf{w})}{(n+1)(1-\beta)^2} v_0.
\end{equation}
At the same time,
\begin{align}\label{estim:ut}
\|\mathbf{a}_{ut}(\mathbf{z})\|^2 &= \left\|2\left(\frac{\|\mathbf{v}\|^2}{n+1}\mathbf{e}  - \mathbf{v}^2\right)  - \frac{v_0}{n+1}\mathbf{e}\right\|^2 \nonumber\\
&= 4\left\| \frac{\|\mathbf{v}\|^2}{n+1}\mathbf{e}  - \mathbf{v}^2\right\|^2+\frac{v_0^2 n}{(n+1)^2}-4\frac{v_0}{n+1}\sum_{i=1}^n \left(\frac{\|\mathbf{v}\|^2}{n+1} - v_i^2\right)\nonumber\\
&\le 4\|\mathbf{v}\|^4 +\frac{v_0^2}{n+1}  + \frac{4 v_0\|\mathbf{v}\|^2}{(n+1)^2}
\le \left(4+\frac{1}{n+1} + \frac{4}{(n+1)^2} \right) v_0^2 \le  \frac{11}{2}v_0^2,
\end{align}
where in the first inequality, we use the bound
$\left\| \frac{\|\mathbf{v}\|^2}{n+1}\mathbf{e}  - \mathbf{v}^2\right\|^2\le\|\mathbf{v}^2\|^2\le\|\mathbf{v}\|^4$.

From the upper bounds (\ref{estim:cr}) and (\ref{estim:ut}), we obtain
\begin{equation}\label{estim:ac}
\frac{1}{2} \| \mathbf{a}_{ac}(\mathbf{z})\|^2 \leq \|\mathbf{a}_{ut}(\mathbf{z})\|^2 + \|\mathbf{a}_{c}(\mathbf{z})\|^2
\le \frac{11}{2}v_0^2 + \frac{\beta^2\rho(\mathbf{w})}{(n+1)(1-\beta)^2} v_0.
\end{equation}
Combining (\ref{estim:alpha_p-a}) and (\ref{estim:ac}), we have
\begin{align*}
1-\alpha &\leq \frac{1-\beta}{\beta \rho(\mathbf{w})}
\left[\|\mathbf{v}\|^2 + \frac{4\kappa}{\sigma^2} 
\left( 11 v_0^2 + \frac{2\beta^2\rho(\mathbf{w})}{(n+1)(1-\beta)^2}v_0\right) \right]\nonumber\\
&\le \frac{1-\beta}{\beta \rho(\mathbf{w})} \left[ \frac{\|\mathbf{v}\|^2}{v_0}+ \frac{44\kappa}{\sigma^2}v_0 \right]v_0 + \frac{8 \kappa }{\sigma^2} \frac{\beta}{(1-\beta)(n+1)} v_0\nonumber\\
&\le \left( \frac{(1-\beta)(n+1)}{\beta (v_0-\|\mathbf{v}\|^2)} \left[ \frac{\|\mathbf{v}\|^2}{v_0}+ \frac{44\kappa}{\sigma^2}v_0 \right] + \frac{8 \kappa }{\sigma^2} \frac{\beta}{(1-\beta)(n+1)}\right) v_0.
\end{align*}

Let us consider the greedy direction in the $\mathbf{w}$-space. Denote
$$ \mathbf{w}^{(0)} = \left(v_0^{(0)},\mathbf{v}^{(0)}\right), \mathbf{w}^{(k)} = \tau_k \mathbf{w}^{(0)}, \; \tau_0=1, \; \tau_{k+1} = (1-\alpha_k)\tau_k, \; k \geq 0.$$
Then $v_0^{(k)} = \tau_k v_0^{(0)}$, and we obtain
\begin{align*}
1-\alpha_k 
&\le \left( \frac{(1-\beta)(n+1)}{\beta (\tau_k v_0^{(0)}-\tau_k^2\|\mathbf{v}^{(0)}\|^2)} \left[ \frac{\tau_k^2\|\mathbf{v}^{(0)}\|^2}{\tau_k v_0^{(0)}}+ \frac{44\kappa}{\sigma^2}\tau_k v_0^{(0)} \right] + \frac{8 \kappa }{\sigma^2} \frac{\beta}{(1-\beta)(n+1)}\right) v_0^{(k)}\nonumber \\
&\le \left( \frac{(1-\beta)(n+1)}{\beta (v_0^{(0)}-v_0^{(k)})} \left[ 1+ \frac{44\kappa}{\sigma^2} v_0^{(0)} \right] + \frac{8 \kappa }{\sigma^2} \frac{\beta}{(1-\beta)(n+1)}\right) v_0^{(k)}.
\end{align*}
Thus, if $v_0$ is small enough, then the method converges quadratically.

\begin{theorem}\label{loc_quadr_conv_ac}
Let $v_0^{(0)} \leq \frac{\sigma^2}{4\kappa}$. Then, for the AC IPA, we have
$$v_0^{(k+1)} \leq \left( \frac{12(1-\beta)(n+1)}{\beta\left(v_0^{(0)} - v_0^{(k)}\right)}  + \frac{2 \beta}{(1-\beta) v_0^{(0)}(n+1)} \right) \left(v_0^{(k)}\right)^2, $$
where $v_0^{(k)}$ is the first coordinate of the vector $\mathbf{w}$ in $k^{th}$ iteration.
\end{theorem}

\section{Numerical results}

Let us present numerical results for randomly generated LCP problems. 
We use a particular version of the generator of random LCP problems from \cite{ParLCP}. For the reader's convenience, we present below its full description.

A necessary and sufficient condition for the solvability of monotone LCP is the existence of a strictly feasible primal-dual pair $(\hat{\mathbf{x}}, \hat {\mathbf{s}})$ such that
$\hat{\mathbf{x}}, \hat{\mathbf{s}} \in \mathbb{R}^n_{+}: \;  - M \hat{\mathbf{x}} + \hat{\mathbf{s}}= \mathbf{q}$.
This condition serves as a starting element for our random generator. It performs the following four steps.
\begin{enumerate}
\item
Choose randomly $\hat{\mathbf{x}}, \hat{\mathbf{s}} \in \mathbb{R}^n_{+}$ with components uniformly distributed in $(0,1)$.
\item
For a random matrix $A \in \mathbb{R}^{n \times n}$ and a lower-triangular matrix $L \in \mathbb{R}^{n \times n}$ with components uniformly distributed in $[0,1]$.
\item
Set $M = A A^T + \eta (L - L^T)$, where $\eta \geq 0$ is a parameter. Define $\mathbf{q} = \hat{\mathbf{s}} - M \hat{\mathbf{x}}$.
\end{enumerate}

Our generator has two parameters $n$ and $\eta$. For our experiments, we choose $\eta = 10$. In Table \ref{tab1}, we present the computational results for different dimensions of three methods: the method from \cite{ParLCP}, justified by the general theory of self-concordant functions (column {\em General IPA}), and two variants of Algorithm \ref{alg:PTS} with the right-hand side (\ref{UTD}) (column {\em UTD IPA}) and right-hand side (\ref{autocorr1}) (column {\em AC IPA}). In the cells of Table \ref{tab1}, we present the average numbers of predictor and corrector steps for the series of twenty-five random problems. In Algorithm 1, we use parameters $\beta = {1 \over 4}$ and $\tau = 1.5$. The problems are solved up to accuracy $\varepsilon = 10^{-7}$.

\begin{table}[h!]
\centering
\caption{Average nr.~of predictor and corrector steps}
\begin{tabular}{|r|r|r|r|}
\hline 
$n$ &  \mbox{General IPA} & \mbox{UTD IPA} & \mbox{AC IPA}\\
\hline
16 & 10.3 / 23.7& 11.2 / 30.7& 8.4 / 20.3 \\
32 & 11.4 / 25.6&12.2 / 31.6 & 9.4 / 21.8 \\
64 & 12.4 / 27.8& 13.8 / 35.8& 11.6 / 32.8\\
128 & 14.6 / 40.4& 17.7 / 46.6& 15.0 / 37.1 \\
256 & 19.0 / 51.2& 21.2 / 68.5& 17.7 / 40.4\\
512 & 22.4 / 67.5& 24.5 / 63.0& 22.7 / 59.0\\
\hline
\end{tabular}\label{tab1}
\end{table}

We can see that the AC IPA is always better than the UTD IPA. Its actual advantage can be seen in the end of the process, where it ensures local quadratic convergence (see Theorem \ref{loc_quadr_conv_ac}). This is confirmed by Table \ref{tab2}, where we present a typical pattern for the progress in the accuracy for all three methods. For instance, at the 23rd predictor iteration the AC IPA improves the accuracy by 3 orders of magnitude. These results correspond to a random problem with $n = 256$ and $\varepsilon = 10^{-8}$ and $\left(\mathbf{x}^{(0)}\right)^T\!\mathbf{s}^{(0)}=66$. As we can see, the first two methods ensure only a local linear rate of convergence.

\begin{table}[h!]
\centering
\caption{Progress in the accuracy} \label{tab2}
\begin{tabular}{|c|c|c|c|}
\hline
\mbox{Accuracy} & \mbox{General IPA} & \mbox{UTD IPA} & \mbox{AC IPA}\\
\hline
\phantom{-1}1 & 12 / 20 & 11 / 19& 11 / 19\\
\phantom{-1}0 & 15 / 26& 15 / 26& 14 / 24\\
\phantom{1}-1 & 17 / 30& 18 / 32& 17 / 30\\
\phantom{1}-2 & 19 / 34& 21 / 38& 19 / 34\\
\phantom{1}-3 & 20 / 36& 23 / 43& 21 / 38\\
\phantom{1}-4 & 21 / 38& 24 / 45& 22 / 40\\
\phantom{1}-5 & 22 / 40& 25 / 47& \\
\phantom{1}-6 & 23 / 42& 26 / 49& \\
\phantom{1}-7 & 24 / 47& 27 / 52& 23 / 43\\
\phantom{1}-8 & 25 / 52& 28 / 55&\\
-12 & & & 24 / 43\\
\hline
\end{tabular}
\text{\vspace{0.2cm}Column Accuracy contains the order of magnitude of $\mathbf{x}^T\mathbf{s}$}{}
\end{table}

\section{Conclusion}

We proposed two parabolic target-space interior-point algorithms for solving monotone linear complementarity problems. The first method is based on a universal tangent direction, while the second one is an auto-correcting version of the previous one.  We proved that both algorithms have the best known worst-case complexity bounds. In the case of the auto-corrector interior-point algorithm, we proved the local quadratic convergence of the method under a non-degeneracy assumption. We also provided numerical results, where we compared the new algorithms with a general method, recently developed for weighted monotone linear complementarity problems.

\section*{Funding}
This research was supported by the National Research, Development and Innovation Office (NKFIH) under grant number  2024-1.2.3-HU-RIZONT-2024-00030. 
The research of Marianna E.-Nagy and Petra Ren\'ata Rig\'o was individually supported by the J\'anos Bolyai Research Scholarship of the Hungarian Academy of Sciences (2024-2027 and 2025-2028, respectively).

\section*{Declarations} 

\textbf{Conflicts of interest} The authors declare that they have no conflict of interest.

\bibliographystyle{abbrv}
\bibliography{MOOR_ref}

\end{document}